\documentclass[11pt]{article}
\usepackage{latexsym,amssymb}
\usepackage{amsmath, amsbsy}
\usepackage{amsopn, amstext}
\usepackage{graphicx}
\newtheorem{lemma}{Lemma}
\newtheorem{theorem}{Theorem}

\newtheorem{proposition}{Proposition}

\newcommand{\diag}{{\text{\upshape diag}}}
\def\NN{\hbox{\rlap{I}\kern.16em N}}
\def\NC{\hbox{\rlap{\kern.24em\raise.1ex\hbox
                  {\vrule height1.3ex width.9pt}}C}}

\title{On the Condition Number of the Total Least Squares
Problem\footnote{The first author was supported by
National Basic Research Program of China  2011CB302400 and National
Science Foundation of China (No. 11071140), and the second
author was supported by Specialized Research
Fund for the Doctoral Program of Higher Education (No. 20070200009).}}
\author{Zhongxiao Jia\\
Department of Mathematical Sciences, Tsinghua University\\
Beijing 100084, P. R. China \\
jiazx@tsinghua.edu.cn\\ \and Bingyu Li\footnote{Corresponding author.}\\
School of Mathematics and Statistics, Northeast Normal University\\
Changchun 130024, P. R. China \\
mathliby@gmail.com  \\
}

\date{}

\begin{document}

\maketitle

\begin{abstract}
This paper concerns singular value
decomposition (SVD)-based computable formulas and bounds for the
condition number of the Total Least Squares (TLS) problem. For the TLS
problem with the coefficient
matrix $A$ and the right-hand side $b$, a new closed formula
is presented for the condition number.
Unlike an important result in the literature that uses the SVDs of both
$A$ and $[A,\ b]$, our formula only requires the SVD of $[A,\ b]$.
Based on the closed formula, both lower and upper bounds for the condition
number are derived. It is proved that they are always sharp and estimate the
condition number accurately. A few lower and
upper bounds are further established that involve at most the smallest
two singular values of $A$ and of $[A,\ b]$. Tightness of these bounds is
discussed, and numerical experiments are presented to confirm our theory and
to demonstrate the improvement of our upper bounds over the two
upper bounds due to Golub and Van Loan as well as Baboulin and Gratton.
Such lower and upper bounds are particularly useful for large scale TLS
problems since they require the computation of only a few singular values
of $A$ and $[A, \ b]$ other than all the singular values of them.

\vskip 5pt \noindent {\bf Keywords:}  total least squares, perturbation,
condition number, singular value decomposition.

\vskip 5pt \noindent {\bf AMS subject classification (2000):} 65F35, 15A12,
15A18.

\end{abstract}

\section{Introduction}

 For given $A \in \mathbb{R}^{m \times n} (m> n)$, $b \in \mathbb{R}^m$,
 the total least squares (TLS) problem can be formulated
as (see, e.g., \cite{Bjorck:1996,GolubVanLoan:1980,PaigeStrakos:2002})
\begin{equation}\label{ScaledTLSprob}
  \min \|[E,\,f]\|_F,
  \text{\quad subject to  \quad}
   b + f \in \mathcal{R}(A+E)  ,
\end{equation}
where $\|\cdot\|_F$ denotes the Frobenius norm of a matrix and
$\mathcal{R}(\cdot)$ denotes the range space. Suppose that $[E_{TLS},
\,f_{TLS}]$ solves the above problem. Then $x = x_{TLS}$ that
satisfies the equation $(A+E_{TLS}) x = b + f_{TLS}$ is called the
TLS solution of (\ref{ScaledTLSprob}). The TLS problem is a formulation of
the linear approximation problem $Ax \approx b$. In this paper,
we concentrate on the inconsistent linear approximation problem, i.e.,
$b \notin{\mathcal{R}(A)}$. Otherwise, $[E_{TLS},\,f_{TLS}]=O$, the zero matrix.

Given a problem, the condition number measures the worst-case
sensitivity of its solution to small perturbations
in the input data. It is well known that the condition number is
independent of perturbations themselves
and is expressed by some information about the original data.
Combined with backward errors, it provides a (possibly approximate)
linear upper bound for the forward error, i.e., the difference
between a perturbed solution and the exact solution.
Since the 1980's, algebraic perturbation analysis for
the TLS problem has received considerable attention;
see \cite{FierroBunch:1996,GolubVanLoan:1980,Liu:1996,WeiM:1992}
and the references therein. From the expressions of perturbation bounds presented
in \cite{FierroBunch:1996,GolubVanLoan:1980,Liu:1996,WeiM:1992},
we can see that there are some essential distinctions between them and
a standard form. The perturbation
bound in \cite{GolubVanLoan:1980} is a standard one in the sense
that it is expressed as some perturbation independent factor
times backward errors. So, this factor
is naturally an upper bound for the TLS condition number. The perturbation
bound in \cite{WeiM:1992} is nonstandard and unusual since the perturbation
bound is not zero when a perturbation is exactly zero. Actually, a careful
observation reveals that the bound is never less than a certain positive
constant under the assumption $b \notin{\mathcal{R}(A)}$.
As a result, it makes no sense to extract
an upper bound for the TLS condition number from this perturbation bound.
The perturbation bounds in \cite{FierroBunch:1996,Liu:1996} are very
different from a standard perturbation bound in that they contain
some information about the perturbed TLS problem, e.g., the TLS solution
$\tilde{x}_{TLS}$ of the perturbed TLS problem in \cite{FierroBunch:1996}
and the right singular vectors associated with the smallest singular values of
both $[A,\ b]$ {\em and} its perturbed matrix in \cite{Liu:1996}. Because of
these features and the fact that the condition number itself has nothing to do
with perturbations, it is impossible to extract upper bounds for the TLS
condition number from the perturbation bounds in \cite{FierroBunch:1996,Liu:1996}.
If one attempts to find suitable upper bounds for the TLS condition number,
some further and complicated treatments are required and it is necessary to
make each of their perturbation bounds become
a standard one, that is, some perturbation independent factor times
the backward error.

In recent years, asymptotic perturbation
analysis and TLS condition numbers have been investigated.
Zhou et al. \cite{ZhouLinWeiQiao:2008} and the authors \cite{LiJia:2009} have
presented a first order perturbation analysis of the TLS problem and established
Kronecker product-based condition number formulas.
Baboulin and Gratton \cite{BaboulinGratton:2010} have derived
an SVD-based closed formula for the TLS condition number, which involves
all the singular values and the right singular vectors of both
$A$ and $[A,\, b]$, and an upper bound, which involves only several
singular values of $A$ and $[A,\, b]$. To our best knowledge, however,
there has been {\em no} lower bound available for the TLS condition number in
the literature.

It is well known that the TLS solution $x_{TLS}$ involves the smallest
singular value and the corresponding right singular
vector of $[A,\, b]$, see, e.g., \cite{GolubVanLoan:1980}.
Very recently, a new classification
has been proposed in \cite{HnetynkovaPlesingerSimaStrakosHuffel:2011} for
the TLS problem in $AX\approx B$ with $B \in \mathbb{R}^{m \times d}$ and
$d\geq 1$. It is based on properties of the SVD of the extended matrix $[B,\ A]$
and has established further results on existence and uniqueness of the TLS
solution. In this paper, based on the intimate relation between SVDs and TLS
problems and motivated by the work of \cite{BaboulinGratton:2010},
we continue our work in \cite{LiJia:2009} to study SVD-based
TLS condition number theory. We will derive a number of results.
Firstly, we establish a new closed formula of the TLS condition number.
It is distinctive that, unlike the result
in \cite{BaboulinGratton:2010} that requires the SVDs of both $A$ and $[A,\ b]$,
our formula only uses the singular values and the right
singular vectors of $[A,\, b]$. Secondly, starting with the
closed formula, we present both {\em lower} and {\em upper}
bounds for the condition number that involve the
singular values of $[A,\, b]$ and the last entries of the right
singular vectors of $[A,\, b]$. Furthermore, we prove that these bounds are
always sharp and can estimate the condition number accurately.
We then focus on cheaply computable
bounds for the TLS condition number. We establish {\em lower} and {\em upper}
bounds that involve {\em at most the smallest two singular values}
of $A$ and $[A,\ b]$. We discuss how tight the bounds are.
These bounds are particularly useful
for large scale TLS problems since they require to compute only very few of the
smallest singular values of $A$ and $[A, \ b]$ rather than
all the singular values of them. So we can compute these bounds by using some
iterative solvers for large SVDs, e.g., \cite{jianiu03,jianiu10}.
From \cite{GolubVanLoan:1980}, as mentioned previously,
an upper bound for the TLS condition number can
be extracted. It has been simplified and applied to evaluate the conditioning
of the TLS problem in \cite{BjorckHeggernesMatstoms:2000}.
We will present numerical experiments to demonstrate
improvements of our upper bounds over the two upper bounds due to
Golub and Van Loan \cite{GolubVanLoan:1980} and Baboulin and Gratton
\cite{BaboulinGratton:2010}, respectively.

We mention that for given $A$ and $b$ the standard least squares (LS) problem
is always and can be much better conditioned than the corresponding TLS problem;
see, e.g., \cite[p.180]{Bjorck:1996}. The results in this paper allow us to
compare the sensitivity of solution of the standard LS problem
to the sensitivity of the solution of the TLS problem.
So it may be better to solve the LS problem if possible.
This is the case when all the errors are confined to the
``observation'' $b$ but $A$ is assumed to be free of errors.
However, this assumption may be unrealistic: sampling errors,
human errors, modeling errors and instrument errors
often imply inaccuracies of $A$ as well. If both $A$ and $b$ are subject to
errors, a reasonable way to
take the errors in $A$ into account may be to introduce perturbations
also in $A$. The TLS problem (\ref{ScaledTLSprob}) is just a natural
formulation for this purpose. We refer the reader to
\cite{Bjorck:1996,GolubVanLoan:1980,HuffelVandewalle:1991}
for more on the introduction of the TLS problem.

The paper is organized as follows. In Section~2, we present some
preliminaries necessary. In Section 3, we establish some useful and
necessary results related to a specific orthogonal matrix.
In Section 4, we present a new closed formula for the TLS condition number.
The bounds for the TLS condition number are derived in Section 5. In Section 6,
we report numerical experiments to show the tightness of our bounds for
the TLS condition number and improvements over Golub-Van
Loan's bound and Baboulin-Gratton's bound. We conclude the paper with
some remarks and future work in Section 7.

Throughout the paper, for given positive integers $m, n$,
denote by $\mathbb{R}^n$ the space of $n$-dimensional real column
vectors, by $\mathbb{R}^{m \times n}$ the space of all $m
\times n$ real matrices, and by $\|\cdot\|$ and $\|\cdot\|_F$ the
2-norm and Frobenius norm of their arguments, respectively.  Given a
matrix $A$, $A(1:i,1:j)$ is a Matlab notation that denotes the submatrix
in the intersection of rows $1, \ldots, i$ and columns $1, \ldots, j$,
and $\sigma_i(A)$ denotes the
$i$th largest singular value of $A$. For a vector $a$, $a(i)$
denotes the $i$th component of $a$, and ${\rm{diag}}(a)$ is a
diagonal matrix whose diagonals are $a(i)$'s.
$I_n$ denotes the $n \times n$ identity matrix and $O_{m,n}$ denotes the
$m \times n$ zero matrix with $O$ a zero matrix whose
dimension is clear from the context. For the matrices $A =[a_1, \ldots, a_n]
= [a_{ij}] \in \mathbb{R}^{m \times n} $ and $B$, $A \otimes B
=[a_{ij} B] $ is the Kronecker product of $A$ and $B$, and the linear
operator ${\rm{vec}}: \mathbb{R}^{m \times n} \rightarrow
\mathbb{R}^{mn} $ is defined by ${\rm{vec}}(A) = [a^T_1, \ldots,
a^T_n]^T$.

\section{Preliminaries}

Throughout the paper, let $A  = \hat{U}
{\rm{diag}}(\hat{\sigma}_1, \ldots, \hat{\sigma}_n)\hat{V}^T$ be the thin SVD
of $A \in \mathbb{R}^{m \times n}$, where  $\hat{\sigma}_1 \geq
\cdots \geq \hat{\sigma}_n$, $\hat{U}\in \mathbb{R}^{m \times n}$,
$\hat{V} \in \mathbb{R}^{n \times n}$ and $\hat{U}^T \hat{U} = I_{n}$,
$\hat{V}^T \hat{V} = I_{n}$, and let $[A,\, b] = U
{\rm{diag}}(\sigma_1, \ldots, \sigma_{n+1}) V^T$ be the thin SVD of $[A,\, b]
\in \mathbb{R}^{m \times (n+1)}$, where $\sigma_1 \geq
\cdots \geq \sigma_{n+1}$, $U = [u_1, \ldots, u_{n+1}]\in
\mathbb{R}^{m \times (n+1)}$,   $V = [v_1, \ldots,
v_{n+1}] \in \mathbb{R}^{(n+1) \times (n+1)}$ and
$U^T U = I_{n+1}$, $V^T V = I_{n+1}$.

The TLS problem (\ref{ScaledTLSprob}) may not have a solution, but it does have
a unique solution if the following condition  holds \cite{PaigeStrakos:2002}:
\begin{equation}\label{PaigeCondition}
A\,\, \text{has rank}\,\, n\,\, \text{and}\,\, b \not\perp \mathcal{U}_{min},
\end{equation}where $\mathcal{U}_{min}$ denotes the left singular vector subspace
of $A$ corresponding to its smallest singular value.
Throughout the paper, we always assume that (\ref{PaigeCondition}) holds.

It is noted in \cite{PaigeStrakos:2002} that condition (\ref{PaigeCondition})
means $\sigma_{n+1} < \hat{\sigma}_{n}$, the existence and
uniqueness condition of the TLS solution given in \cite{GolubVanLoan:1980}.
Under the condition that $\sigma_{n+1} < \hat{\sigma}_{n}$, it is proved in
\cite{GolubVanLoan:1980} that \begin{eqnarray}\label{STLSsolutionExp}
 x_{TLS} & = & (A^T A - \sigma^2_{n+1} I )^{-1} A^T b \\ \label{SVDSTLSsol}
 &=& -
\left[\frac{v_{n+1}(1)}{v_{n+1}(n+1)}, \ldots,
\frac{v_{n+1}(n)}{v_{n+1}(n+1)} \right]^T.
\end{eqnarray}
We comment that (\ref{PaigeCondition}) implies that $A^T b \neq 0$.
So $x_{TLS} \neq 0$.

Given the TLS
problem (\ref{ScaledTLSprob}), let $\tilde{A}= A + \Delta A$,
$\tilde{b} = b+ \Delta b$, where $\Delta A$ and $\Delta b$ denote
the perturbations in $A$ and $b$, respectively. Consider the
perturbed TLS problem
\begin{equation}\label{perturbedScaledTLS}
\min \|[E,\,f]\|_F \,\,\,{\text{subject to}} \,\, \tilde{b}+f \in
\mathcal{R}(\tilde{A} + E).
\end{equation}
Under the assumption that $b \notin \mathcal{R}(A)$, it follows from
(\ref{PaigeCondition}) that $0 < \sigma_{n+1} < \hat{\sigma}_{n}$.
In \cite{LiJia:2009}, the following result is established
for the TLS solution $\tilde{x}_{TLS}$ of the perturbed TLS problem
(\ref{perturbedScaledTLS}).

\begin{theorem}\label{MyProp}
Suppose that the TLS problem {\rm(\ref{ScaledTLSprob})} satisfies
$0 < \sigma_{n+1} < \hat{\sigma}_{n}$. Define $r = A x_{TLS} - b$ and
$G = [x^T_{TLS}, \,\,-1] \otimes I_m \in \mathbb{R}^{m \times (mn+m)}$. If
$\|[\Delta A, \,\, \Delta b]\|_F$ is small enough, then the perturbed
problem {\rm (\ref{perturbedScaledTLS})} has a unique TLS solution
$\tilde{x}_{TLS}$. Moreover,
\begin{equation}\label{vecExpression}
\tilde{x}_{TLS} = x_{TLS} + K~ \left[
\begin{array}{c}
{\rm{vec}}(\Delta A) \\
\Delta b \\
\end{array}
\right]
+ \mathcal{O}(\|[\Delta A, \,\, \Delta b]\|^2_F),
\end{equation}
where \begin{equation}\label{bigK} K = \left(A^T A - \sigma^2_{n+1}
I_{n} \right)^{-1} \left( 2 A^T \frac{r}{\|r\|} \frac{r^T}{\|r\|}
G-A^T G-[I_n \otimes r^T,\, O_{n,m}] \right)\in \mathbb{R}^{n \times (mn +m)}.
\end{equation}
\end{theorem}

It is shown in \cite{GolubVanLoan:1980} that
\begin{equation}\label{optimalfunction}
\sigma^2_{n+1} = \frac{\|r\|^2}{1 + \|x_{TLS}\|^2}
\end{equation}and
\begin{equation}\label{gradienteqn}
A^T r = \frac{\|r\|^2}{1 + \|x_{TLS}\|^2} x_{TLS} = \sigma^2_{n+1} x_{TLS}.
\end{equation}
From (\ref{SVDSTLSsol}), it
follows that
\begin{equation}\label{solusv}
v_{n+1} = \frac{1}{ \sqrt{1 +\|x_{TLS}\|^2}} \left[
                                                 \begin{array}{c}
                                                    x_{TLS} \\
                                                   -1 \\
                                                 \end{array}
                                               \right]
\end{equation}
up to a sign $\pm 1$. We will use the above two relations later.
The following basic properties of the Kronecker products of matrices
can be found in \cite{Graham:1981} and are needed later:
\begin{eqnarray}\nonumber
&& (A_1 \otimes A_3)(A_2\otimes A_4) = ( A_1 A_2) \otimes (A_3 A_4), \\
\nonumber &&
(A_1 \otimes A_2)^T = A^T_1 \otimes A^T_2,
\end{eqnarray}where $A_i, i=1, \ldots, 4$ are matrices of appropriate
sizes.

\section{Some results related to a specific orthogonal matrix}

In this section, we establish a number of results that are related to
a specific orthogonal matrix. They play a central role in deriving
our lower and upper bounds for the condition number of the TLS problem
in Section 5.

\begin{proposition}\label{GolubCSdecomp}
Let $W$ be an arbitrary $(n+1) \times (n+1)$ orthogonal matrix
with $W(n+1, n+1) = - \alpha$, $0< \alpha < 1$.
Denote $W_{11} = W(1:n, 1:n)$. Then
\begin{equation}\label{W11singularvalues}
\sigma_1(W_{11}) = \cdots = \sigma_{n-1}(W_{11}) = 1,\,\,
\sigma_n({W_{11}}) = \alpha.
\end{equation}Furthermore,
$W$ can be written as
\begin{equation}\label{CSdecomp}
W = \left[
            \begin{array}{cc}
              W_{11} & \sqrt{1 - \alpha^2}~ \bar{u}_n \\
              \sqrt{1 - \alpha^2}~ \bar{v}^T_n & - \alpha \\
            \end{array}
          \right],
\end{equation}
where $\bar{u}_n$ and $\bar{v}_n$ are the left and right singular
vectors associated with the smallest singular value of $W_{11}$.
\end{proposition}
{\bf{Proof.}} It is an immediate result of Theorem 2.6.3 in
\cite{GolubVanLoan:1996}. \hfill
$\Box$

Let $[\bar{\beta}_1, \ldots, \bar{\beta}_n, -\alpha]$ be
the last row of $W$.
From (\ref{CSdecomp}) we have
\begin{equation}\label{lastrowV}
\bar{v}^T_n = \frac{1}{\sqrt{1 - \alpha^2}} [\bar{\beta}_1, \ldots,
\bar{\beta}_n].
\end{equation}

Since ($\alpha^{-1}$, $\bar{u}_n$, $\bar{v}_n$) is the
largest singular triplet of $W^{-T}_{11}$, from (\ref{W11singularvalues})
we get the SVD of $W^{-T}_{11}$:
$$W^{-T}_{11} = \alpha^{-1} \bar{u}_n \bar{v}^T_n + \sum^{n-1}_{i=1}
\bar{u}_i \bar{v}^T_i,$$
where $\bar{u}_i, \bar{v}_i, i =1, \ldots, n-1$ represent
the left and right singular vectors associated with the singular value one.
Then, by (\ref{lastrowV}) we obtain
\begin{eqnarray}\label{Nowweuseit}
W^{-T}_{11}&=& \left[\alpha^{-1} \bar{u}_n, \bar{u}_1, \ldots,
\bar{u}_{n-1} \right] \left[
\begin{array}{ccc}
\frac{\bar{\beta}_1}{\sqrt{1 - \alpha^2}} & \cdots  &
\frac{\bar{\beta}_n}{\sqrt{1 - \alpha^2}} \\
\bar{v}_1(1) & \ldots & \bar{v}_1(n) \\
\vdots &  \cdots &  \vdots \\
\bar{v}_{n-1}(1) & \cdots & \bar{v}_{n-1}(n)\\
\end{array}
\right] \\
\label{inV11struct} &=& \left[ \frac{\alpha^{-1} \bar{\beta}_1}{\sqrt{1 -
\alpha^2}}\bar{u}_n + w_1, \ldots, 
\frac{\alpha^{-1} \bar{\beta}_n}{\sqrt{1 - \alpha^2}}\bar{u}_n + w_n
\right],
\end{eqnarray}where $\bar{v}_i(k)$ denotes the $k$th component
of $\bar{v}_i$ and
\begin{equation}
w_k = \sum^{n-1}_{i =1} \bar{v}_i(k)\bar{u}_i,\ k
= 1, \ldots, n. \label{wk}
\end{equation}
From (\ref{Nowweuseit}) we get
\begin{equation}\label{sumbeta}
\sqrt{\sum_{i=1}^{n-1}\bar{v}_i^2(n)}=\frac{\sqrt{1-\alpha^2-\bar{\beta}_n^2}}
{\sqrt{1-\alpha^2}},
\end{equation}which will be used later.

Before proceeding, we need the following lemma.

\begin{lemma}\label{mylovelemma}
For given matrices $A_1, A_2 \in \mathbb{R}^{n \times n}$, if $A^T_1
A_2 = O$, then
\begin{equation}\label{mylovelemmaeqn}
\frac{1}{2} (\|A_1\| + \|A_2\|) \leq \|A_1 + A_2\| \leq \|A_1\| + \|A_2\|.
\end{equation}
\end{lemma}
{\bf{Proof.}} The upper bound in (\ref{mylovelemmaeqn}) is obvious.
It suffices to prove the lower one.
For an arbitrary vector $x \in \mathbb{R}^{n}$, from
$(A_1 x)^T (A_2 x) = 0$ it follows that
\begin{equation}\nonumber
\|A_1 x\|, \|A_2 x\| \leq \|A_1 x + A_2 x\|
\end{equation}and that
\begin{eqnarray}\nonumber
\|A_1\| &=& {\rm{max}}_{\|x\| =1} \|A_1 x\| \leq {\rm{max}}_{\|x\|
=1} \|A_1 x + A_2 x\| = \|A_1 + A_2\|,
\\\nonumber
\|A_2\| &=& {\rm{max}}_{\|x\| =1} \|A_2 x\| \leq {\rm{max}}_{\|x\|
=1} \|A_1 x + A_2 x\| = \|A_1 + A_2\|.
\end{eqnarray}
So, the assertion is proved. \hfill $\Box$

Now we are in a position to show the following two propositions.

\begin{proposition}\label{bigsurpriseprop}
Let $W$ be an arbitrary $(n+1) \times (n+1)$ orthogonal matrix with
$W(n+1, n+1) = - \alpha$, $0< \alpha < 1$.
Let $W_{11} = W(1:n, 1:n)$ and $[\bar{\beta}_1, \ldots,
\bar{\beta}_n, -\alpha]$ be the last row of $W$.
Then for $\bar{S} = {\rm{diag}}(\bar{s}_1, \ldots,
\bar{s}_n)$ with $\bar{s}_1, \ldots, \bar{s}_n$ arbitrary
positive numbers ordered as $0< \bar{s}_1 \leq \bar{s}_2 \leq
\cdots \leq \bar{s}_n$, we have
\begin{eqnarray}
&& \underline{c} := \frac{1}{2} \left( \frac{\alpha^{-1}
\sqrt{\bar{\beta}^2_1 \bar{s}^2_1 + \cdots + \bar{\beta}^2_n \bar{s}^2_n }
}{\sqrt{1 - \alpha^{2}}} + \frac{\sqrt{1 - \alpha^2 - \bar{\beta}^2_n  }}{
\sqrt{1 - \alpha^2} }  \bar{s}_n \right) \nonumber \\
&\leq &
\left\|W^{-T}_{11} \bar{S} \right\|\leq \bar{c} := \frac{\alpha^{-1}
\sqrt{\bar{\beta}^2_1 \bar{s}^2_1 + \cdots + \bar{\beta}^2_n \bar{s}^2_n }
}{\sqrt{1 - \alpha^{2}}} + \bar{s}_n. \label{bigsurprisepropineq}
\end{eqnarray}
\end{proposition}
{\bf{Proof.}} Following (\ref{inV11struct}), we get
\begin{eqnarray}\nonumber 
W^{-T}_{11} \bar{S} = \left[ \frac{\alpha^{-1}
\bar{\beta}_1 \bar{s}_1}{\sqrt{1 - \alpha^2}}\bar{u}_n + \bar{s}_1 w_1,
\ldots,  \frac{\alpha^{-1}
\bar{\beta}_n \bar{s}_n}{\sqrt{1 - \alpha^2}}\bar{u}_n + \bar{s}_n w_n
\right].
\end{eqnarray}
Define the matrices
\begin{equation}\nonumber
A_1 = \left[\frac{\alpha^{-1} \bar{\beta}_1 \bar{s}_1}{\sqrt{1 -
\alpha^2}}\bar{u}_n, \ldots, \frac{\alpha^{-1} \bar{\beta}_n
\bar{s}_n}{\sqrt{1 - \alpha^2}}\bar{u}_n \right]
\end{equation}
and
\begin{equation}\label{mydef}
A_2 =W^{-T}_{11}\bar{S} - A_{1} =
\left[\bar{s}_1 w_1, \ldots, \bar{s}_n w_n   \right].
\end{equation}
Then
\begin{equation}\label{newsum}
\|A_1\| = \frac{\alpha^{-1}}{\sqrt{1 - \alpha^2}} \sqrt{\bar{\beta}^2_1
\bar{s}^2_1+ \cdots+ \bar{\beta}^2_n \bar{s}^2_n}.
\end{equation}
From (\ref{wk}) and (\ref{sumbeta}) we obtain
\begin{equation}\nonumber
\|w_n\| = \sqrt{\sum^{n-1}_{i=1} \bar{v}^2_i(n)}=
\frac{\sqrt{1-\alpha^2-\bar{\beta}^2_n}}{\sqrt{1-\alpha^2}}\,\mbox{ and }\,
\|[w_1, \ldots, w_n]\| =\left\|\sum^{n-1}_{i=1} \bar{u}_i \bar{v}^T_i \right\|=1.
\end{equation}
Furthermore, since $\|\bar{S}\| = \bar{s}_n$, it follows from
\begin{equation}\nonumber
\|\bar{s}_n w_n\| \leq \|A_2\| \leq \left\|[w_1, \ldots, w_n]
\right\| \|\bar{S}\|
\end{equation}
that
\begin{equation}\label{surprise3}
\frac{\sqrt{1-\alpha^2 - {\bar{\beta}}^2_n}}{\sqrt{1-\alpha^2}}\bar{s}_n
\leq \|A_2\| \leq \bar{s}_n.
\end{equation}
Note that
$A^T_1 A_2 = O$.
Based on Lemma \ref{mylovelemma} and combining  (\ref{mydef}), (\ref{newsum}) with
(\ref{surprise3}), we establish the desired inequality. \hfill
$\Box$

\begin{proposition}\label{oneremark}
Suppose that $0<\alpha \leq\frac{1}{2}$, where $\alpha$ is defined as in
Proposition \ref{bigsurpriseprop}. Then for $\underline{c}$ and
$\bar{c}$ in Proposition \ref{bigsurpriseprop}, we have
\begin{equation}\label{remarkineqn}
\underline{c} < \bar{c}< 4 \underline{c}.
\end{equation}
\end{proposition}
{\bf{Proof.}} If
$\frac{|\bar{\beta}_n|}{\sqrt{1-\alpha^2}}<\frac{\sqrt{3}}{2}$, it is
easy to verify that
\begin{equation}\nonumber
\frac{\sqrt{1-\alpha^2 - \bar{\beta}^2_n}}{\sqrt{1-\alpha^2}}
>\frac{1}{2}
\end{equation}
and that
\begin{equation}\nonumber
        \underline{c} >\frac{1}{4} \bar{c}.
\end{equation}
Thus, (\ref{remarkineqn}) holds. If
$\frac{|\bar{\beta}_n|}{\sqrt{1-\alpha^2}}\geq\frac{\sqrt{3}}{2}$, then
$$
\alpha^{-1} \frac{|\bar{\beta}_n|}{\sqrt{1-\alpha^2}}
\geq\frac{\sqrt{3}}{2}\alpha^{-1} >1,
$$
so $\alpha^{-1} \frac{|\bar{\beta}_n|}{\sqrt{1-\alpha^2}} \bar{s}_n >
\bar{s}_n$, from which and the definitions of $\bar{c}$ and
$\underline{c}$ it follows that
\begin{eqnarray}\nonumber
\bar{c} &< & \frac{\alpha^{-1} \sqrt{\bar{\beta}^2_1 \bar{s}^2_1 + \cdots
+ \bar{\beta}^2_n \bar{s}^2_n }  }{\sqrt{1 - \alpha^{2}}} + \alpha^{-1}
\frac{|\bar{\beta}_n|}{\sqrt{1-\alpha^2}} \bar{s}_n
\\\nonumber
&\leq & \frac{2 \alpha^{-1} \sqrt{\bar{\beta}^2_1 \bar{s}^2_1 + \cdots +
\bar{\beta}^2_n \bar{s}^2_n }  }{\sqrt{1 - \alpha^{2}}}
\\\nonumber
&\leq & \frac{2 \alpha^{-1} \sqrt{\bar{\beta}^2_1 \bar{s}^2_1 + \cdots +
\bar{\beta}^2_n \bar{s}^2_n }  }{\sqrt{1 - \alpha^{2}}} +  \frac{2 \sqrt{1
- \alpha^2 - \bar{\beta}^2_n  }}{ \sqrt{1 - \alpha^2} }  \bar{s}_n = 4
\underline{c}.
\end{eqnarray}
Thus, (\ref{remarkineqn}) holds. \hfill $\Box$

This proposition means that the upper bound is at most four times of the lower
bound in Proposition~\ref{bigsurpriseprop}. So we can estimate
$\left\|W^{-T}_{11} \bar{S} \right\|$ accurately by its lower or upper bound
within no more than four times of the exact $\left\|W^{-T}_{11} \bar{S} \right\|$.

\section{A closed formula for the TLS condition number}

Throughout the paper, we follow the definition of condition number
in \cite{GohbergKoltracht:1993,Rice:1966}. Let $g: \mathbb{R}^p
\longrightarrow \mathbb{R}^q$ be a continuous map in normed linear
spaces defined on an open set $D_g \subset \mathbb{R}^p$. For a
given $a_0 \in D_g$, $a_0 \neq 0$, with $g(a_0) \neq 0$, if $g$ is
differentiable at $a_0$, then the absolute condition number of $g$ at $a_0$ is
\begin{equation}\label{Rice_abs}\kappa_g(a_0) = \|g'(a_0)\|,
\end{equation} and the relative condition number is
\begin{equation}\label{Rice_rel}
\kappa^{\text{rel}}_g(a_0) =  \frac{\|g'(a_0)\| \|a_0\|}{\|g(a_0)\|},
\end{equation}
where $g'(a_0)\in \mathbb{R}^{q \times p}$ denotes the Jacobian of $g$ at $a_0$.

In \cite{BaboulinGratton:2010}, an SVD-based closed formula for the condition
number of the TLS problem was presented. Denote by $\kappa_{g}(A,b)$ the
absolute TLS condition number. It was shown in
\cite{BaboulinGratton:2010} that
\begin{equation}\label{Baboulin}
\kappa_{g}(A,b) =\sqrt{ 1+\|x_{TLS}\|^2} \left\| \hat{D}~[\hat{V}^T,
\,O_{n,1}]~  V~ [D, \,O_{n,1}]^T \right\|,
\end{equation}where
\begin{eqnarray}\nonumber
\hat{D} &= & {\rm{diag}}\left(\frac{1}{\hat{\sigma}^2_1 -
\sigma^2_{n+1} }, \ldots, \frac{1}{\hat{\sigma}^2_n - \sigma^2_{n+1}
} \right), \\ \nonumber D &=& {\rm{diag}} \left(
\sqrt{\sigma^2_1 + \sigma^2_{n+1}}, \ldots, \sqrt{\sigma^2_n +
\sigma^2_{n+1}}\right).
\end{eqnarray}

Next we will derive a new SVD-based formula for the TLS condition
number. It is distinctive that, unlike (\ref{Baboulin}) that involves the
singular values and right singular vectors of both $[A,\,b]$ and $A$,
our formula only uses those of $[A,\,b]$.

Denote $a
={\rm{vec}}(A)$ and define the following function in a small
neighborhood of $[a^T,b^T]^T \in \mathbb{R}^{m (n+1)}$:
\begin{equation}\nonumber
\begin{array}{ccc}
  g: \mathbb{R}^{m (n+1)}  &  \longrightarrow & \mathbb{R}^n \\
  \small{\left[
    \begin{array}{c}
      \tilde{a} \\
      \tilde{b} \\
    \end{array}
  \right]}
   & \longmapsto & \tilde{x}_{TLS}= (\tilde{A}^T \tilde{A} -
   \tilde{\sigma}^2_{n+1}I_n)^{-1} \tilde{A}^T
  \tilde{b},
\end{array}
\end{equation}
where $\tilde{A} = A + \Delta A$, $\tilde{a} = {\rm{vec}}(\tilde{A}) =
a + {\rm{vec}}(\Delta A)
$, $\tilde{b} = b+ \Delta b$, and $\tilde{\sigma}_{n+1} = \sigma_{n+1}
([\tilde{A},\,\tilde{b}])$.
Then we have
$g([a^T, b^T]^T) = x_{TLS}$. 
Based on Theorem \ref{MyProp}, we can present the following result.

\begin{theorem}\label{knoneckercondform}
Given the TLS problem {\rm (\ref{ScaledTLSprob})}, let $\kappa_{g}(A,b)$
and $\kappa^{\text{rel}}_{g}(A,b)$ be the absolute and
relative condition numbers of the TLS problem, respectively. Then
\begin{equation}\label{knoneckercond}
\kappa_{g}(A,b) = \|K\|, \,\,\kappa^{\rm rel}_{g}(A,b) = \frac{\|K\| \|[A,
\,b]\|_F}{\|x_{TLS}\|},
\end{equation}
where $K$ is defined as in {\rm (\ref{bigK})}.
\end{theorem}
{\bf{{Proof}}}.\,\,
Recall that our TLS problem satisfies $0 < \sigma_{n+1} < \hat{\sigma}_n$.
By Theorem \ref{MyProp} and the definition of $g$, we see that $g$ is
differentiable at $[a^T,b^T]^T$ and $g'\left([a^T,b^T]^T \right) = K$. Then the
assertion follows from (\ref{Rice_abs}) and (\ref{Rice_rel}) . \hfill
$\Box$

The formulas for the TLS problem in Theorem \ref{knoneckercondform}
depend on Kronecker products of matrices. We comment that the formula for
$\kappa^{\rm rel}_{g}(A,b)$ has the same form as that for
$\kappa^K_{STLS}$ with $\lambda =1$ in Theorem 3.3 of \cite{LiJia:2009}, and
as stated in \cite{LiJia:2009}, mathematically we have $\kappa^{\rm rel}_{g}(A,b)=
\kappa^M_{STLS}$ when $\lambda =1$ in Theorem 3.1 of \cite{ZhouLinWeiQiao:2008},
where $\kappa^M_{STLS}$ is the relative condition number of the Scaled TLS problem
with $\lambda$ the scaling factor that is derived in \cite{ZhouLinWeiQiao:2008}.

Now we can establish our computable formula for the TLS
condition number.

\begin{theorem}\label{myformationforcond}
Given the TLS problem {\rm (\ref{ScaledTLSprob})}, let $[A,\,b]  =
U {\rm{diag}}(\sigma_1, \ldots, \sigma_{n+1}) V^T$
be the thin SVD of $[A,\,b]$ with $V_{11} = V(1:n, 1:n)$. Then
\begin{equation}\label{mynewclosed}
\kappa_{g}(A,b) =\sqrt{ 1+\|x_{TLS}\|^2} ~\|V^{-T}_{11} S\|,
\end{equation}
where $S = {\rm{diag}}(s_1, \ldots, s_n)$ with $s_i =
\frac{\sqrt{\sigma^2_i + \sigma^2_{n+1}}}{\sigma^2_i -
\sigma^2_{n+1}}$, $i = 1, \ldots, n$.
\end{theorem}

{\bf{Proof.}}
Consider expression (\ref{bigK}) of $K$. By the
properties of Kronecker product of matrices, we get

\begin{equation}\nonumber
G  G^T = \left( [ x^T_{TLS}, \,-1] \otimes I_m  \right)  \left(
\left[
                                \begin{array}{c}
                                  x_{TLS} \\
                                  -1 \\
                                \end{array}
                              \right] \otimes I_m  \right)
= (1+\|x_{TLS}\|^2) I_m,
\end{equation}
\begin{equation}\nonumber
[I_n \otimes r^T, \,O_{n,m}]  G^T = [ I_n \otimes r^T, \,O_{n,m}] \left[
                            \begin{array}{c}
                               x_{TLS} \otimes I_m \\
                              - I_m \\
                            \end{array}
                          \right]
= (I_n \otimes r^T) ( x_{TLS} \otimes I_m) =
 x_{TLS}\, r^T
\end{equation}
and
\begin{equation}\nonumber
[I_n \otimes r^T, \,O_{n,m}]  \left[
                                      \begin{array}{c}
                                        I_n \otimes r \\
                                        O_{m,n} \\
                                      \end{array}
                                    \right]
 =(I_n \otimes r^T) (I_n \otimes r) = \|r\|^2 I_n.
\end{equation}
Thus, we have
\begin{eqnarray}\nonumber
&& \left( 2 A^T \frac{r}{\|r\|} \frac{r^T}{\|r\|}G  - A^T G -
[I_n \otimes r^T, \,O_{n,m}]  \right) \\ \nonumber && \cdot \left( 2
G^T\frac{r}{\|r\|} \frac{r^T}{\|r\|} A -  G^T A -
 \left[\begin{array}{c}
 I_n \otimes r \\
 O_{m,n} \\
 \end{array}
 \right] \right) \\\nonumber
&=& (1+\|x_{TLS}\|^2) A^T A + \|r\|^2 I_n - x_{TLS} r^T A - A^T r x^T_{TLS} \\
\nonumber &=&  ( 1+\|x_{TLS}\|^2  ) A^T A + \|r\|^2 I_n- 2 \sigma^2_{n+1}
x_{TLS} x^T_{TLS},
 \end{eqnarray}
where the last equality uses the relation $A^T r x^{T}_{TLS}=
\sigma^2_{n+1} x_{TLS} x^T_{TLS}$, which is obtained from (\ref{gradienteqn}).
Denote $P = A^T A - \sigma^2_{n+1}I_n$. We get
\begin{eqnarray}\nonumber
 K  K^T &=& P^{-1} \left( (  1+\|x_{TLS}\|^2 ) A^T A +
\|r\|^2 I_n- 2 \sigma^2_{n+1} x_{TLS} x^T_{TLS}   \right) P^{-1}
\\ \label{firstchange}
&=& (1+\|x_{TLS}\|^2) P^{-1} \left(  A^T A + \sigma^2_{n+1} I_n - \frac{2
\sigma^2_{n+1} x_{TLS} x^T_{TLS}}{1+\|x_{TLS}\|^2} \right) P^{-1}\\ \label{forcond}
&=& (1+\|x_{TLS}\|^2) \left(P^{-1} + 2 \sigma^2_{n+1} P^{-1} \left( I_n  -
 \frac{ x_{TLS} x^T_{TLS}}{1 + \|x_{TLS}\|^2} \right) P^{-1} \right),
\end{eqnarray}
where the second equality used (\ref{optimalfunction}).
Denote $V=[v_1, \ldots, v_{n+1}]$ and note that
\begin{eqnarray}\nonumber
[A,\,b]^T [A,\,b] - \sigma^2_{n+1} I_{n+1}& =& \sum^{n+1}_{i=1}
\sigma^2_i v_i v^T_i - \sigma^2_{n+1} \sum^{n+1}_{i=1} v_i v^T_i \\
\nonumber &=& \sum^{n}_{i=1} (\sigma^2_i -\sigma^2_{n+1} )v_i v^T_i.
\end{eqnarray}
We get
\begin{eqnarray}\nonumber
P &=& [I_n, O] \sum^{n}_{i=1}
(\sigma^2_i -\sigma^2_{n+1} )v_i v^T_i \left[
\begin{array}{c}
I_n \\
O \\
\end{array}
\right] \\ \nonumber &=& [I_n, O] [v_1, \ldots, v_n]
\left[
                                     \begin{array}{ccc}
                                       \sigma^2_1 - \sigma^2_{n+1} &  &  \\
                                        & \ddots &  \\
                                        &  & \sigma^2_n - \sigma^2_{n+1} \\
                                     \end{array}
                                   \right]
\left[
                                             \begin{array}{c}
                                               v^T_1 \\
                                               \vdots \\
                                               v^T_n \\
                                             \end{array}
                                           \right]
 \left[
                                                      \begin{array}{c}
                                                        I_n \\
                                                        O \\
                                                      \end{array}
                                                    \right] \\ \label{forP}
&=& V_{11}  \diag(\sigma^2_1 - \sigma^2_{n+1}, \ldots, \sigma^2_n - \sigma^2_{n+1})
V^T_{11} := V_{11} \Lambda V^T_{11}.
                                   \label{pdecom}
\end{eqnarray}
Similarly,   since $v_{n+1} = \frac{1}{\sqrt{1 +
\|x_{TLS}\|^2}} \left[
                                          \begin{array}{c}
                                            x_{TLS} \\
                                            -1 \\
                                          \end{array}
                                        \right]
$ (c.f. (\ref{solusv})), we have
\begin{eqnarray}\nonumber
 I_{n+1} - \frac{1}{{1 + \|x_{TLS}\|^2}} \left[
 \begin{array}{cc}
 x_{TLS} x^T_{TLS} & -x_{TLS} \\
 -x_{TLS} & 1 \\
 \end{array}
 \right] = I_{n+1} - v_{n+1} v^T_{n+1} =[v_1, \ldots, v_n] \left[
 \begin{array}{c}
 v^T_1 \\
 \vdots \\
 v^T_n \\
 \end{array}
 \right]
\end{eqnarray}and
\begin{equation}\label{forV11}
 I_n - \frac{x_{TLS} x^T_{TLS}}{1 + \|x_{TLS}\|^2 } = V_{11} V^T_{11}.
\end{equation}
Combining
(\ref{forP}) and (\ref{forV11}), we have
\begin{eqnarray}\nonumber
&& P^{-1} + 2 \sigma^2_{n+1} P^{-1} \left( I_n - \frac{ x_{TLS} x^T_{TLS}}{1 +
\|x_{TLS}\|^2} \right) P^{-1} \\ \nonumber
&=& V^{-T}_{11} \Lambda^{-1}
V^{-1}_{11} + 2 \sigma^2_{n+1} \left( V^{-T}_{11} \Lambda^{-1}
V^{-1}_{11} \right)  V_{11} V^T_{11} \left( V^{-T}_{11}
\Lambda^{-1}V^{-1}_{11} \right) \\ \label{Wewilluseit}
&=&
V^{-T}_{11} \Lambda^{-1} V^{-1}_{11} + 2 \sigma^2_{n+1} V^{-T}_{11}
\Lambda^{-2} V^{-1}_{11} \\ \nonumber  
&=& V^{-T}_{11} \left( \Lambda^{-1} + 2
\sigma^2_{n+1} \Lambda^{-2} \right) V^{-1}_{11} = \left( V^{-T}_{11}
S \right) \left( V^{-T}_{11} S \right)^T.
\end{eqnarray}
Then it follows from (\ref{forcond}) and Theorem \ref{knoneckercondform}
that the desired equality holds.
 \hfill
$\Box$

By Theorem \ref{myformationforcond}, we can calculate $\kappa_{g}(A,b)$ by
solving a linear system with the coefficient matrix $V^T_{11}$.
Next we show that the condition number of
$V^T_{11}$ is exactly $\sqrt{1 + \|x_{TLS}\|^2}$.

\begin{theorem}\label{svforV11}
Under the conditions of Theorem \ref{myformationforcond}, we have
\begin{equation}\label{V11singularvalues}
\sigma_1(V_{11}) = \cdots = \sigma_{n-1}(V_{11}) = 1,
\sigma_n({V_{11}}) = \frac{1}{\sqrt{1 + \|x_{TLS}\|^2}}
\end{equation}
and
\begin{equation}\label{V11conditionnumber}
\kappa(V_{11}) = \frac{\sigma_1(V_{11})}{\sigma_n(V_{11})}=
\sqrt{1 + \|x_{TLS}\|^2}.
\end{equation}
\end{theorem}
{\bf{Proof.}} By (\ref{solusv}) we get $V(n+1,
n+1)=-\frac{1}{\sqrt{1 + \|x_{TLS}\|^2}}$. Recalling that
$x_{TLS} \neq 0$, we get $\frac{1}{\sqrt{1 + \|x_{TLS}\|^2}} < 1 $.
Applying Proposition \ref{GolubCSdecomp}, we prove the theorem.
 \hfill $\Box$

\section{Bounds for the TLS condition number}

\subsection{Sharp lower and upper bounds based on the SVD of $[A,\,b]$ }

In this subsection, we apply Theorem
\ref{myformationforcond} and show how to estimate
$\kappa_{g}(A,b)$ accurately without computing
$\left\| V^{-T}_{11} S \right\|$.

For the TLS problem (\ref{ScaledTLSprob}), from now on we denote
$\alpha = - V(n+1, n+1)$. By (\ref{solusv}) we get
\begin{equation}\label{alpha}
\alpha = \frac{1}{\sqrt{1+\|x_{TLS}\|^2}}.
\end{equation}
Recalling that $x_{TLS} \not=0$, we have $0 < \alpha < 1$.
Then we have the following theorem.

\begin{theorem}\label{theorem5}
For the TLS problem {\rm(\ref{ScaledTLSprob})}, it holds that
\begin{equation}\label{notmyloveestimate}
\underline{\kappa}:= \alpha^{-1} \frac{\sqrt{\sigma^2_n + \sigma^2_{n+1}}}
{\sigma^2_n -\sigma^2_{n+1}}  \leq \kappa_{g}(A,b) \leq
\bar{\kappa}:=  \alpha^{-2} \frac{\sqrt{\sigma^2_n + \sigma^2_{n+1}}}
{\sigma^2_n -\sigma^2_{n+1}}.
\end{equation}
\end{theorem}
{\bf{Proof.}} As before, let $[A,\,b]=
U {\rm{diag}}(\sigma_1, \ldots, \sigma_{n+1}) V^T$
be its thin SVD with $V_{11} = V(1:n,1:n)$.
 From
(\ref{V11singularvalues}) and (\ref{alpha}) it follows that
\begin{equation}
\sigma_n(V^{-T}_{11}) = \left(\sigma_1(V_{11}) \right)^{-1} =1,\,\,
\sigma_1(V^{-T}_{11}) = \left(\sigma_n(V_{11}) \right)^{-1} = \alpha^{-1}.
\end{equation}
Define $S = {\rm{diag}}(s_1, \ldots, s_n)$ with $s_i =
\frac{\sqrt{\sigma^2_i + \sigma^2_{n+1}}}{\sigma^2_i -
\sigma^2_{n+1}}$, $i = 1, \ldots, n$.
 We then have
 \begin{equation}\nonumber
  s_n  = \sigma_n(V^{-T}_{11}) \|S\| \leq\|V^{-T}_{11} S\| \leq
  \|V^{-T}_{11}\| \|S\| = \alpha^{-1} s_n.
 \end{equation}
Therefore, by Theorem \ref{myformationforcond} we get the
desired inequality. \hfill $\Box$

We see that in (\ref{notmyloveestimate}) the ratio of the upper bound and the
lower bound is $\frac{1}{\alpha}$. As a consequence,
both the bounds estimate $\kappa_g(A,b)$ within $\frac{1}{\alpha}$ times.
Therefore, if $\alpha\in (0,1)$ is not small, say $\frac{1}{2}<\alpha<1$, both
the bounds are very tight and they estimate $\kappa_g(A,b)$ accurately.

Starting with (\ref{Baboulin}), Baboulin and
Gratton \cite[Corollary 1]{BaboulinGratton:2010} have derived the following
upper bound
\begin{equation}\label{bgbound}
\kappa_g(A,b)\leq \alpha^{-1}\frac{\sqrt{\sigma_1^2+\sigma_{n+1}^2}}
{\hat{\sigma}_n^2-\sigma_{n+1}^2}:=\bar{\kappa}(A,b),
\end{equation}
which uses $\hat{\sigma}_n-\sigma_{n+1}$ to estimate the conditioning of
$x_{TLS}$. The smaller  $\hat{\sigma}_n-\sigma_{n+1}$ is, the possibly worse
conditioned the TLS problem is. It has two distinctions with our lower
and upper bounds in (\ref{notmyloveestimate}). First, (\ref{bgbound}) involves
the SVDs of both $A$ and $[A,\ b]$ while (\ref{notmyloveestimate}) only makes
use of that of $[A,\ b]$. Second, since $ \hat{\sigma}^2_n -\sigma^2_{n+1}
\leq \sigma^2_n - \sigma^2_{n+1}$ and $\sigma_1\geq\sigma_n$, our lower and
upper bounds can be considerably more accurate than (\ref{bgbound}) for $\alpha$
not small. We now present a family of examples to illustrate it and the tightness of
the bounds in (\ref{notmyloveestimate}) for $\frac{1}{2}<\alpha<1$.

\begin{table}[h]
\begin{center}
{ \doublerulesep18.0pt \tabcolsep0.01in
\begin{tabular}{ |c| c| c | c |c |c|c |}
\hline
$e_p$ &$\alpha$   &
$\hat{\sigma}_n - \sigma_{n+1} $& $\kappa_g(A,b)$ &
$~\alpha^{-1} \frac{\sqrt{\sigma^2_n + \sigma^2_{n+1}}}{\sigma^2_n -
\sigma^2_{n+1}}~$ &
 $~\alpha^{-2} \frac{\sqrt{\sigma^2_n + \sigma^2_{n+1}}}{\sigma^2_n -
\sigma^2_{n+1}} ~$ &~$\alpha^{-1}\frac{\sqrt{\sigma_1^2+\sigma_{n+1}^2}}
{\hat{\sigma}_n^2-\sigma_{n+1}^2}~$
 \\
\hline
$10^{-3}$ & $~0.516~$ &
 $9.89\times 10^{-4}$ & $1.38\times 10^3$
 &$1.37\times 10^{3}$
 &$2.66\times 10^{3}$
 &$4.16\times 10^5$
 \\
 \hline
$10^{-7}$ &  $~0.792~$   &
 $9.99\times 10^{-8}$ & $8.93\times 10^6$
  &$8.93\times 10^{6}$
 &$1.13\times 10^{7}$
 &$2.54\times 10^9$
 \\
 \hline
$10^{-10}$ &  $~0.859~$   &
 $1.00\times 10^{-10}$ & $ 8.24\times 10^{9}$
  &$8.23\times 10^{9}$
 &$9.59\times 10^{9}$
 &$2.33\times 10^{12}$
 \\
 \hline
\end{tabular}}
\end{center}
\end{table}

{\bf{Example.}}
We construct TLS problems as in \cite[Example 1]{BaboulinGratton:2010}:
Define
$$
[A, \,b] = Q \left[
                 \begin{array}{c}
                   \Sigma \\
                   O \\
                 \end{array}
               \right]
 V^T \in \mathbb{R}^{m \times (n+1)}, \,Q = I_m - 2 y y^T, \,
 V = I_{n+1} - 2 z z^T,
 $$ where $y\in \mathbb{R}^{m}$ and $z \in \mathbb{R}^{n+1}$ are random
 unit vectors, and $\Sigma = \diag(n, n-1, \ldots, 1, 1-e_p)$ for a given
 parameter $e_p$. Note that $e_p=\sigma_n-\sigma_{n+1}$. We have $e_p\geq
 \hat{\sigma}_n-\sigma_{n+1}$. By taking small values of $e_p$,
 we get different TLS problems whose conditioning becomes worse and condition
 number becomes larger as $e_p$ becomes smaller. Fixing $m = 100$, $n=20$, and
 taking $e_p = 10^{-3}, 10^{-7}, 10^{-10}$, respectively,
 we get three different TLS problems whose solutions are computed by the
 SVD of $[A,\ b]$ and (\ref{SVDSTLSsol}). As indicated by the results of
 $\kappa_g(A,b)$, as $e_p$ decreases, the TLS problem becomes worse conditioned.
 This is also reflected by the decay of $\hat{\sigma}_n - \sigma_{n+1}$;
 see (\ref{bgbound}), and Theorems~\ref{myfirstestimate}--\ref{improveTLSthm}
 and \cite{GolubVanLoan:1980}. Since the $\alpha$'s are
 bigger than $0.5$ and not small, the lower and upper bounds in
 (\ref{notmyloveestimate}) estimate the TLS condition numbers accurately,
 and they are much sharper than bound (\ref{bgbound}) by roughly two to three
 orders.

In view of (\ref{notmyloveestimate}) and the comments after its proof
as well as the above example, it is only possible and significant to
improve the bounds essentially for the case that $\alpha$ is small relative to
one. Without loss of generality, we assume that
\begin{equation}\label{alphaassumption}
0<\alpha\leq\frac{1}{2}.
\end{equation}
It will appear that we can establish some lower bound $\underline{\kappa}$ and
upper bound $\bar{\kappa}$ such that
$\underline{\kappa}<\bar{\kappa}<4\underline{\kappa}$
holds. As a result, together with (\ref{notmyloveestimate}),
we can estimate the TLS condition number $\kappa_g(A,b)$ accurately.

\begin{theorem}\label{bigsurprisethm}
Given the TLS problem {\rm (\ref{ScaledTLSprob})}, let $[A,\,b]=
U {\rm{diag}}(\sigma_1, \ldots, \sigma_{n+1}) V^T$
be its thin SVD. Denote $[\beta_1, \ldots,
\beta_n, -\alpha]$ be the last row of $V$ and
$S = {\rm{diag}} (s_1, \ldots, s_n)$, $s_i =
\frac{\sqrt{\sigma^2_i + \sigma^2_{n+1}}}{\sigma^2_i -
\sigma^2_{n+1}}$, $i = 1, \ldots, n$. Then
\begin{eqnarray}\nonumber 
&& \underline{\kappa} := \frac{1}{2} \left( \frac{\alpha^{-2}
\sqrt{\beta^2_1 s^2_1 + \cdots + \beta^2_n s^2_n }  }{\sqrt{1 -
\alpha^{2}}} + \frac{\sqrt{1 - \alpha^2 - \beta^2_n  }}{ \sqrt{1 -
\alpha^2} } \alpha^{-1} s_n \right)
\\\nonumber
&\leq & \kappa_{g}(A,b)\leq \bar{\kappa} := \frac{\alpha^{-2}
\sqrt{\beta^2_1 s^2_1 + \cdots + \beta^2_n s^2_n }  } {\sqrt{1 -
\alpha^{2}}} + \alpha^{-1}s_n.
\end{eqnarray}Moreover, if $0< \alpha \leq \frac{1}{2}$, then
\begin{equation}\nonumber
\underline{\kappa} < \bar{\kappa} < 4 \underline{\kappa}.
\end{equation}
\end{theorem}
{\bf{Proof.}} Recall that $0 < \alpha < 1$.
Noticing that $0< s_1 \leq s_2 \leq \cdots \leq s_n$
and applying Proposition \ref{bigsurpriseprop}, we have
\begin{eqnarray}\nonumber
&& \frac{1}{2} \left( \frac{\alpha^{-1} \sqrt{\beta^2_1 s^2_1 + \cdots + \beta^2_n s^2_n }  }{\sqrt{1 - \alpha^{2}}} + \frac{\sqrt{1
- \alpha^2 - \beta^2_n  }}{ \sqrt{1 - \alpha^2} } s_n \right) \\
\nonumber 
&\leq & \left\|V^{-T}_{11} S \right\| \leq \frac{\alpha^{-1}
\sqrt{\beta^2_1 s^2_1 + \cdots + \beta^2_n s^2_n } }{\sqrt{1 -
\alpha^{2}}} + s_n,
\end{eqnarray}where $V_{11} = V(1:n,1:n)$.
By Theorem \ref{myformationforcond}, we get the first part of the
theorem. Furthermore, we obtain the second part of the theorem  by
Proposition \ref{oneremark}. \hfill $\Box$

The significance of this theorem is that we can estimate the
condition number $\kappa_g(A,b)$ accurately by its lower or upper
bound without calculating $\|V^{-T}_{11}S\|$, i.e., solving the matrix
equation $V_{11}^TW=S$ for $W$ and computing the 2-norm of $W$,
which is expensive when $n$ is large.

\subsection{Lower and upper bounds based on a few singular values of
$A$ and $[A,\,b]$}

In \cite{Malyshev:2003}, some bounds for the condition number of the
Tikhonov regularization solution have been established using {\em only} a few
singular values of $A$, where $A$ is the coefficient matrix of the
least squares problem under consideration.
Such kind of bounds are particularly appealing
for large scale TLS problems, because the condition number in
Theorem \ref{myformationforcond} and the bounds
in Theorem~\ref{bigsurprisethm} involve all the singular values of $[A, b]$
and are impractical for computational purpose.

Actually, we have presented such kind of bound (\ref{notmyloveestimate}),
but as we commented there, for small $\alpha$,
the bounds may overestimate or underestimate $\kappa_g(A,b)$.
In the following, we establish some new lower and upper bounds in the
same spirit and finally achieve sharper lower and upper bounds
for the case of $0<\alpha\leq\frac{1}{2}$.

\begin{theorem}\label{firstestimateforcond}
For the TLS problem {\rm (\ref{ScaledTLSprob})}, we have
\begin{equation}\label{luestimate}
\underline{\kappa}_1 \leq \kappa_{g}(A,b) \leq  \bar{\kappa}_1,
\end{equation}
where
\begin{eqnarray}\label{lubound}
\underline{\kappa}_1 = \alpha^{-1}\frac{\sqrt{
\hat{\sigma}^2_{n-1} + \sigma^2_{n+1} }}{ \hat{\sigma}^2_{n-1} -
\sigma^2_{n+1}}, \,\, \bar{\kappa}_1 = \alpha^{-1}\frac{\sqrt{\hat{\sigma}^2_{n}
+ \sigma^2_{n+1} }} {  \hat{\sigma}^2_{n}
- \sigma^2_{n+1}}.
\end{eqnarray}
\end{theorem}

{\bf{Proof.}} Recall that $r= A x_{TLS} -b$.
Denote
$$M = (A^T A - \sigma^2_{n+1} I_n)^{-1} \left( (
1+\|x_{TLS}\|^2) A^T A +  \|r\|^2 I_n \right) (A^T A - \sigma^2_{n+1}
I_n)^{-1}.
$$
From (\ref{firstchange}), we have
\begin{equation}\label{fromotherproof}
K  K^T = M - 2 \sigma^2_{n+1} (A^T A - \sigma^2_{n+1} I_n)^{-1} x_{TLS}
x^T_{TLS} (A^T A - \sigma^2_{n+1} I_n)^{-1}.
\end{equation}
From now on, denote by $\lambda_i(M)$ the $i$th algebraically largest
eigenvalue of $M$, where $M$ is an arbitrary symmetric matrix. By
the Courant-Fischer theorem \cite[p. 182]{Horn:1999}, we get
\begin{equation}\label{inequality1}
 \lambda_2(M) \leq  \lambda_1(K K^T).
\end{equation}
Furthermore, since $2\sigma^2_{n+1}(A^T A - \sigma^2_{n+1} I_n)^{-1}
x_{TLS} x^T_{TLS} (A^T A - \sigma^2_{n+1} I_n)^{-1}$ is nonnegative definite,
the following inequality holds
\begin{equation}\label{inequlity2}
 \lambda_1(K  K^T) \leq \lambda_1(M).
\end{equation}
Combining (\ref{inequality1}) with (\ref{inequlity2}) and based on
(\ref{knoneckercond}), we have
\begin{equation}\nonumber
\sqrt{\lambda_2(M)}  \leq \kappa_{g}(A,b) \leq \sqrt{\lambda_1(M)}.
\end{equation}

It is easy to verify that the eigenvalues of $M$ form the set
\begin{equation}\nonumber
\left\{ \frac{(1+\|x_{TLS}\|^2) \hat{\sigma}^2_j +
 \|r\|^2}{(\hat{\sigma}^2_j - \sigma^2_{n+1})^2}
\right\}^{n}_{j=1}.
\end{equation}
We define the function
\begin{equation}\nonumber
h(\sigma) = \frac{( 1+\|x_{TLS}\|^2) {\sigma}^2 +
 \|r\|^2}{({\sigma}^2 - \sigma^2_{n+1})^2}, \,\,\sigma >
\sigma_{n+1}
\end{equation}
and differentiate it to get
\begin{equation}\nonumber
h'(\sigma) = \frac{-2 \sigma^3 (1+\|x_{TLS}\|^2) - 2 \sigma ( 1+\|x_{TLS}\|^2) \sigma^2_{n+1} - 4 \sigma \|r\|^2}{(\sigma^2 -
\sigma^2_{n+1})^3}.
\end{equation}
It is seen that $h'(\sigma) < 0$ and $h(\sigma)$ is decreasing in
the interval $(\sigma_{n+1}, \infty) $. Thus, we get
\begin{equation}\nonumber
 \lambda_1(M) = \frac{(1+\|x_{TLS}\|^2)
\hat{\sigma}^2_n + \|r\|^2}{(\hat{\sigma}^2_n -
\sigma^2_{n+1})^2},\,\, \lambda_2(M) = \frac{(1+\|x_{TLS}\|^2)
\hat{\sigma}^2_{n-1} + \|r\|^2}{(\hat{\sigma}^2_{n-1} -
\sigma^2_{n+1})^2}
\end{equation}
and
\begin{equation}\nonumber
\frac{\sqrt{(1+\|x_{TLS}\|^2) \hat{\sigma}^2_{n-1} +
\|r\|^2}}{\hat{\sigma}^2_{n-1} - \sigma^2_{n+1}} \leq
\kappa_{g}(A,b) \leq
    \frac{\sqrt{( 1+\|x_{TLS}\|^2) \hat{\sigma}^2_n +
\|r\|^2}}{\hat{\sigma}^2_n - \sigma^2_{n+1}}.
\end{equation}
Recalling (\ref{alpha}) and that $\frac{\|r\|^2}{1 + \|x_{TLS}\|^2} =
\sigma^2_{n+1}$ completes the proof. \hfill $\Box$

{\bf Remark}.
Since $\hat{\sigma}_n\leq \sigma_n$ and $\sigma_n\leq\sigma_1$, we have
proved that
\begin{equation}\label{compare}
\bar{\kappa}_1 \leq \alpha^{-1}\frac{\sqrt{\sigma_1^2+\sigma_{n+1}^2}}
{\hat{\sigma}_n^2-\sigma_{n+1}^2}:=\bar{\kappa}(A,b),
\end{equation}
which is just the upper bound (\ref{bgbound}) derived by Baboulin and
Gratton \cite{BaboulinGratton:2010}. Therefore, our upper bound
$\bar{\kappa}_1$ in (\ref{lubound}) is always sharper than
Baboulin-Gratton's bound. Moreover, the improvement must be
significant when $\frac{\sigma_1}{\hat{\sigma}_n}>1$ considerably.

It is seen that the lower and upper bounds in
Theorem \ref{firstestimateforcond} are marginally different provided
that  $\hat{\sigma}_n$ and $\hat{\sigma}_{n-1}$ are close. This
means that in this case both bounds are very tight. For the case
that $\hat{\sigma}_n$ and $\hat{\sigma}_{n-1}$ are not close, we
next give a new lower bound that can be better than that in Theorem
\ref{firstestimateforcond}.

\begin{theorem}\label{myfirstestimate}
It holds that
\begin{equation}\nonumber
\underline{\kappa}_2 \leq \kappa_{g}(A,b)\leq \bar{\kappa}_1,
\end{equation}where $\bar{\kappa}_1$ is defined as in
Theorem \ref{firstestimateforcond} and
\begin{equation}\nonumber
\underline{\kappa}_2 =\frac{1}{\alpha\sqrt{
{\hat{\sigma}}^2_n
 - \sigma^2_{n+1}}}.
\end{equation}Moreover, when $\hat{\sigma}_{n-1} \geq \sigma_{n+1} +
\sqrt{{\hat{\sigma}}^2_n - \sigma^2_{n+1}}$, we have
\begin{equation}\nonumber
\underline{\kappa}_1 \leq \underline{\kappa}_2. \end{equation}
\end{theorem}

{\bf{Proof.}} Denote $P = A^T A - \sigma^2_{n+1}I_n$. From
(\ref{forcond}), we have
\begin{equation}\nonumber
\frac{1}{1+\|x_{TLS}\|^2}
 ~K  K^T =P^{-1} + 2\sigma^2_{n+1} P^{-1}
 \left( I_n - \frac{x_{TLS} x^T_{TLS}}{1+\|x_{TLS}\|^2} \right) P^{-1}.
\end{equation}
Since the second term in the right-hand side of the above
relation is positive definite, we have
\begin{equation}\nonumber
(1+\|x_{TLS}\|^2)\lambda_1(P^{-1}) \leq \lambda_1(  ~K K^T ),
\end{equation}
that is,
\begin{equation}\nonumber
\frac{1+\|x_{TLS}\|^2}{{\hat{\sigma}^2_n} - \sigma^2_{n+1}}\leq
\kappa^2_{g}(A,b),
\end{equation}
which used (\ref{knoneckercond}).
Thus, recalling (\ref{alpha}), we obtain the first part of the theorem.

The second part of the theorem is obtained by noting that
\begin{equation}\nonumber
\frac{ \sqrt{ \hat{\sigma}^2_{n-1} + \sigma^2_{n+1} }}{
\hat{\sigma}^2_{n-1} - \sigma^2_{n+1}}< \frac{1}{ \hat{\sigma}_{n-1}
- \sigma_{n+1}}\leq \frac{1} {\sqrt{{\hat{\sigma}}^2_n -
\sigma^2_{n+1}}}
\end{equation}under the assumption that $\hat{\sigma}_{n-1} -\sigma_{n+1}\geq
\sqrt{{\hat{\sigma}}^2_n - \sigma^2_{n+1}}$. \hfill $\Box$

{\bf{Remark 1.}} At a first glance, the assumption in the second part
of the theorem seems not so direct but we can justify that it indeed
implies that $\hat{\sigma}_n$ and $\hat{\sigma}_{n-1}$ are not
close. Actually, it is direct to verify that the second part of Theorem
\ref{myfirstestimate} holds under the slightly stronger but much
simpler condition that
\begin{equation}\nonumber
 \hat{\sigma}_{n-1} \geq 2 \hat{\sigma}_n.
 \end{equation}

{\bf{Remark 2.}} From
\begin{equation}\nonumber
\frac{\bar{\kappa}_1}{\underline{\kappa}_2} =
\frac{\sqrt{\hat{\sigma}^2_n +  \sigma^2_{n+1}}}
{\sqrt{\hat{\sigma}^2_n - \sigma^2_{n+1}}} = \sqrt{\frac{1 +
\frac{\sigma^2_{n+1}}{\hat{\sigma}^2_n}} {1 -  \frac{\sigma^2_{n+1}}
{\hat{\sigma}^2_n}  }},
\end{equation}it is seen that $\frac{\bar{\kappa}_1}{\underline{\kappa}_2}>
1$  provided that $\sigma_{n+1} > 0$. Only for $\sigma_{n+1} = 0$,
$\bar{\kappa}_1 = \underline{\kappa}_2$ holds. Then, $b \in
\mathcal{R}(A)$ and $r = 0$.

We observe from  the above Remark 2 that the bounds for $\kappa_g(A,b)$
in Theorem~\ref{myfirstestimate} are tight when
$\frac{\sigma_{n+1}}{\hat{\sigma}_n}$ is considerably smaller than one.
On the other hand, if $\frac{\sigma_{n+1}}{\hat{\sigma}_n} $ is
not small, these bounds may not be tight. For this case, we will
present a new upper bound for better estimating $\kappa_g(A,b)$.

Keep {\rm (\ref{alpha})} and (\ref{alphaassumption}) in mind. Based on
Propositions \ref{bigsurpriseprop}--\ref{oneremark}, we establish the
following theorem.

\begin{theorem}\label{improveTLSthm}
If $0 < \alpha \leq\frac{1}{2}$, then
\begin{equation}\label{TLSimprove}
\underline{\kappa}_2: =\frac{1}{\alpha\sqrt{\hat{\sigma}^2_n - \sigma^2_{n+1}}}
\leq \kappa_g(A,b) < \bar{\kappa}_2:=\sqrt{\frac{1+31 \rho^2}{1-\rho^2}}
 \frac{1}{\alpha\sqrt{\hat{\sigma}^2_n - \sigma^2_{n+1}}},
\end{equation}
where $\rho =\frac{\sigma_{n+1}}{\sigma_n}$.
\end{theorem}

{\bf Proof.}\  The lower bound is the same as that in
Theorem \ref{myfirstestimate}. We only need to prove the right-hand side
of (\ref{TLSimprove}). As before, let $[A,\,b]=
U {\rm{diag}}(\sigma_1, \ldots, \sigma_{n+1}) V^T$
be its thin SVD with $V_{11} = V(1:n,1:n)$.
From (\ref{forcond}), (\ref{pdecom}) and (\ref{Wewilluseit}),
we get
\begin{eqnarray}\nonumber
\frac{1}{1+\|x_{TLS}\|^2}~K K^T &=& P^{-1} + 2 \sigma^2_{n+1} P^{-1}
 \left( I_n - \frac{x_{TLS} x^T_{TLS}}{1+\|x_{TLS}\|^2} \right) P^{-1},
 \\ \label{fornewbound}
 &=& V^{-T}_{11} \Lambda^{-1} V^{-1}_{11} + 2 \sigma^2_{n+1} V^{-T}_{11}
 \Lambda^{-2} V^{-1}_{11}:=P^{-1} + C,
\end{eqnarray}
where $P=A^T A - \sigma^2_{n+1} I_n$, $ \Lambda = {\rm{diag}}
(\sigma^2_1 - \sigma^2_{n+1}, \ldots, \sigma^2_n -
\sigma^2_{n+1})$. Denote
\begin{eqnarray}\nonumber
D &=& {\rm{diag}}(d_1, \ldots, d_n),\, d_i = \frac{\sigma_{n+1}}
{\sigma^2_i - \sigma^2_{n+1}}, i = 1, \ldots, n,
\\\nonumber
T &=& {\rm{diag}}(t_1, \ldots, t_n),\, t_i =
\frac{1}{\sqrt{\sigma^2_i - \sigma^2_{n+1}}}, i = 1, \ldots, n.
\end{eqnarray}
Then $P^{-1} = \left( V^{-T}_{11} T\right) \left( T V^{-1}_{11}
\right)$ and $C = 2 \left( V^{-T}_{11} D\right)  \left(D V^{-1}_{11}
\right)$.

Note that $0 < d_1 \leq d_2 \leq \cdots \leq d_n$ and $0 < t_1 \leq
t_2 \leq \cdots \leq t_n$. From Proposition
\ref{bigsurpriseprop}, we get
\begin{eqnarray}\nonumber 
&& \frac{1}{2} \left( \frac{\alpha^{-1} \sqrt{\beta^2_1 d^2_1 + \cdots +
\beta^2_n d^2_n }  }{\sqrt{1 - \alpha^{2}}} + \frac{\sqrt{1
- \alpha^2 - \beta^2_n  }}{ \sqrt{1 - \alpha^2} } d_n \right) \\
\label{lambdabigsurprise3} &\leq & \left\|V^{-T}_{11} D  \right\|
\leq \frac{\alpha^{-1} \sqrt{\beta^2_1 d^2_1 + \cdots + \beta^2_n
d^2_n }  }{\sqrt{1 - \alpha^{2}}} + d_n
\end{eqnarray}
and
\begin{eqnarray}\nonumber
&& \frac{1}{2} \left( \frac{\alpha^{-1} \sqrt{\beta^2_1 t^2_1 + \cdots + \beta^2_n t^2_n }  }{\sqrt{1 - \alpha^{2}}} + \frac{\sqrt{1
- \alpha^2 - \beta^2_n  }}{ \sqrt{1 - \alpha^2} } t_n \right) \\
\nonumber
 &\leq & \left\|V^{-T}_{11} T  \right\|
\leq \frac{\alpha^{-1} \sqrt{\beta^2_1 t^2_1 + \cdots + \beta^2_n
t^2_n }  }{\sqrt{1 - \alpha^{2}}} + t_n,
\end{eqnarray}
respectively, where  $[\beta_1, \ldots,
\beta_n, -\alpha]$ is the last row of $V$ as defined previously. Define
$k_n = \frac{d_n}{t_n}=\frac{\sigma_{n+1}}{\sqrt{\sigma^2_n -
\sigma^2_{n+1}   }}$. Then
\begin{equation}\nonumber
\frac{d_1}{t_1}= \frac{\sigma_{n+1}}{\sqrt{\sigma^2_1 -
\sigma^2_{n+1}   }} \leq k_n\,, \ldots, \frac{d_{n-1}}{t_{n-1}}
=\frac{\sigma_{n+1}}{\sqrt{\sigma^2_{n-1} - \sigma^2_{n+1}   }} \leq
k_n.
\end{equation}
Thus, from (\ref{lambdabigsurprise3}) we get
\begin{eqnarray}\label{firstimproineqn}
\frac{1}{\sqrt{2}} \|C\|^{\frac{1}{2}} = \left\|V^{-T}_{11} D
\right\|   \leq  k_n \left( \frac{\alpha^{-1}
 \sqrt{\beta^2_1 t^2_1 + \cdots + \beta^2_n t^2_n }  }{\sqrt{1 - \alpha^{2}}} + t_n
 \right).
\end{eqnarray}
For the lower and upper bounds on $\left\|V^{-T}_{11} T
\right\|$ above, Propositions~\ref{bigsurpriseprop}--\ref{oneremark} tell us
that
\begin{eqnarray}\nonumber
\frac{\alpha^{-1} \sqrt{\beta^2_1 t^2_1 + \cdots + \beta^2_n t^2_n }
}{\sqrt{1 - \alpha^{2}}} + t_n &<& 2 \left( \frac{\alpha^{-1}
\sqrt{\beta^2_1 t^2_1 + \cdots + \beta^2_n t^2_n } }{\sqrt{1 -
\alpha^{2}}} + \frac{\sqrt{1 - \alpha^2 - \beta^2_n  }}{ \sqrt{1 -
\alpha^2} } t_n \right) \\ \label{lambdabigsurprise4} &<& 4
\left\|V^{-T}_{11} T \right\|.
\end{eqnarray}
Therefore, based on (\ref{firstimproineqn}) and (\ref{lambdabigsurprise4}), we
obtain
\begin{equation}\nonumber
\frac{1}{\sqrt{2}} \|C\|^{\frac{1}{2}} < 4 k_n \left\|V^{-T}_{11} T
\right\| = 4 k_n \|P^{-1}\|^{\frac{1}{2}}
\end{equation}
and
\begin{equation}\label{improveE}
\|C\| < 32 k^2_n \|P^{-1}\|.
\end{equation}
Combining   (\ref{fornewbound}) with (\ref{improveE}) and based on
(\ref{knoneckercond}), we get
\begin{eqnarray}\nonumber
\kappa_g(A,b) = \|K\| = \|K K^T\|^{\frac{1}{2}} & < & \sqrt{1 +
32k^2_n}\sqrt{1+\|x_{TLS}\|^2} \|P^{-1}\|^{\frac{1}{2}}\\\nonumber
&=& \sqrt{\frac{1+31 \rho^2}{1-\rho^2}}\frac{\sqrt{1+\|x_{TLS}\|^2}}
{\sqrt{\hat{\sigma}^2_n-\sigma^2_{n+1}}}\\\nonumber
&=& \sqrt{\frac{1+31 \rho^2}{1-\rho^2}}\frac{1}
{\alpha\sqrt{\hat{\sigma}^2_n-\sigma^2_{n+1}}},
\end{eqnarray}where the last equality uses (\ref{alpha}).
\hfill $\Box$

{\bf{Remark.}} It
is clear that the bounds in Theorem
\ref{improveTLSthm} are tight when $\rho=\frac{\sigma_{n+1}}{\sigma_{n}}$
is considerably smaller than one.  Note
that the lower and upper bounds in Theorem \ref{myfirstestimate}
differ considerably when $\frac{\sigma_{n+1}}{\hat{\sigma}_n} $ is
close to one.
The result in this theorem is of particular importance when
$\frac{\sigma_{n+1}}{\hat{\sigma}_n} $ is close to one, since the upper bound
here can be considerably sharper than the upper bound of
Theorem \ref{myfirstestimate} when $\rho=\frac{\sigma_{n+1}}{\sigma_{n}}$ is
small.

The improvement of $\bar{\kappa}_2$ over
$\bar{\kappa}_1$ can be illustrated as follows.
For $\rho = \frac{\sigma_{n+1}}{\sigma_{n}}$ small, we have
\begin{eqnarray*}
\bar{\kappa}^{\text{rel}}_2 := \frac{\bar{\kappa}_2}{\|x_{TLS}\|} \|[A,\,b]\|_F &=&
\sqrt{\frac{1+31 \rho^2}{1-\rho^2}}\frac{\sqrt{1+ \|x_{TLS}\|^2}}{\|x_{TLS}\|}
\frac{\|[A,\,b]\|_F}{\sqrt{\hat{\sigma}^2_n - \sigma^2_{n+1} }},
\end{eqnarray*}
an upper bound for $\kappa^{\text{rel}}_g(A,b)$,
is a moderate multiple of $\frac{1}
{\sqrt{\hat{\sigma}^2_n - \sigma^2_{n+1} }}$,
while
\begin{eqnarray*}
\bar{\kappa}^{\text{rel}}_1 := \frac{\bar{\kappa}_1}{\|x_{TLS}\|} \|[A,\,b]\|_F &=&
\frac{\sqrt{1+ \|x_{TLS}\|^2}}{\|x_{TLS}\|} \frac{ \sqrt{\hat{\sigma}^2_n +
\sigma^2_{n+1} }    }{ \hat{\sigma}^2_n - \sigma^2_{n+1} }
\|[A,\,b]\|_F
\end{eqnarray*}
is a moderate multiple of $\frac{1} {\hat{\sigma}^2_n -
\sigma^2_{n+1} }$. So the improvement of $\bar{\kappa}^{\text{rel}}_2$ over
$\bar{\kappa}^{\text{rel}}_1$ becomes significant when ${\sigma_{n+1}}$ and
${\hat{\sigma}_n} $ are close.

Golub and Van Loan \cite{GolubVanLoan:1980} derive an upper bound for the relative
condition number of the TLS problem,
which, in our notation and case, is simplified as
\begin{equation}\label{golubloanbound}
\kappa^{\rm rel}_{TLS}(A,b):=
\frac{9\sigma_1}{\sigma_n-\sigma_{n+1}}\left(1+\frac{\|b\|}
{\hat{\sigma}_n-\sigma_{n+1}}\right)\frac{\|[A,\,b]\|_F}{\|b\|-\sigma_{n+1}}.
\end{equation}
From (\ref{bgbound}),  Babcoulin and Gratton \cite{BaboulinGratton:2010} get the
following upper bound  for the relative condition number for the TLS problem:
\begin{equation}\label{bgrelbound}
\bar{\kappa}^{\rm rel}(A,b):=\frac{\sqrt{1+\|x_{TLS}\|^2}}{\|x_{TLS}\|}
\frac{\sqrt{\sigma_1^2+\sigma_{n+1}^2}}
{\hat{\sigma}_n^2-\sigma_{n+1}^2}\|[A,\ b]\|_F.
\end{equation}

We will numerically illustrate improvements of our bounds
$\bar{\kappa}^{\text{rel}}_1$ and $\bar{\kappa}^{\text{rel}}_2$ over
(\ref{golubloanbound}) and (\ref{bgrelbound}) in the next section.

\section{Numerical experiments}

We present numerical experiments to illustrate the tightness of the bounds
in Theorems \ref{myfirstestimate}--\ref{improveTLSthm} and to show that our upper
bounds can be much better than (\ref{golubloanbound}) and (\ref{bgrelbound}).
For a given TLS problem, the TLS solution is computed by (\ref{SVDSTLSsol}).
All the experiments were run using Matlab 7.8.0 with the machine precision
$\epsilon_{\rm mach}=2.22\times 10^{-16}$ under the Microsoft Windows XP operating
system. Keep $0<\alpha< 1$ in mind. As we have seen from Theorem~\ref{theorem5}
and the comments after it as well as the numerical example there,
for $\alpha$ not small, e.g., $\alpha\in (\frac{1}{2},1)$, we can
estimate the TLS condition number accurately since the lower and
upper bounds in (\ref{notmyloveestimate}) are both
sharp in this case. Next we will be concerned with only the case that $\alpha$
is not near one. For all the test problems, we always have $0<\alpha<\frac{1}{2}$.

{\bf{Example 1.}} The data $A \in \mathbb{R}^{m \times (m-2)}$,
$b \in \mathbb{R}^{m}$ are taken from \cite{HuffelVandewalle:1991}:
\begin{equation*}
A = \left[
      \begin{array}{cccc}
        m-1 & -1 & \cdots & -1 \\
        -1 & m-1 & \cdots & -1 \\
        \vdots & \vdots & \cdots & \vdots \\
        -1 & -1 & \cdots & m-1 \\
        -1 & -1 & \cdots & -1 \\
        -1 & -1 & \cdots & -1 \\
      \end{array}
    \right], \,\,
b= \left[
     \begin{array}{c}
       -1 \\
       -1\\
       \vdots\\
       -1 \\
       m-1\\
       -1 \\
     \end{array}
   \right].
\end{equation*}
So the exact \begin{equation}\nonumber
x_{TLS} = [-1, -1, \ldots, -1]^T \in \mathbb{R}^{m-2},\,\, \hat{\sigma}_n =
\sqrt{2m}, \,\, \sigma_{n+1} = \sqrt{m},\, \,\alpha=\frac{1}{\sqrt{m-1}}.
\end{equation}

\begin{table}[h]
\begin{center}
{ \doublerulesep18.0pt \tabcolsep0.01in
\begin{tabular}{|c| c | c| c | c  |c|c|}
\hline $m$   &
$\kappa^{\text{rel}}_g(A,b)$ &
 $\underline{\kappa}^{\text{rel}}_{2}$&
 $\bar{\kappa}^{\text{rel}}_{2} $  & $\bar{\kappa}^{\text{rel}}_{1} $  & $\bar{\kappa}^{\text{rel}}(A,b)$ &
 $\kappa^{\rm rel}_{TLS}(A,b) $
 \\
\hline
 $~200~$   &
  $2.01\times 10^2$   & $2.00\times 10^2$&
 $2.15\times 10^2$ &$3.46\times 10^{2}$& $~ 2.83\times 10^{3}$ &
 $5.15\times 10^{3}$
 \\
 \hline
  $~500~$   &
  $5.01\times 10^2$
 & $5.00\times 10^2$& $5.15\times 10^2$ &$8.65\times 10^{2}$ & $~ 1.12\times 10^{4}$
 &$1.21\times 10^{4}$
 \\
 \hline
  $~1000~$   &
  $ 1.00\times 10^3$   & $1.00\times 10^3$
 & $1.02\times 10^3$ &$1.73\times 10^{3}$ & $~ 3.16\times 10^{4}$
 &$2.35\times 10^{4}$
 \\
 \hline
\end{tabular}}
\caption{Example 1}\label{tab5}
\end{center}
\end{table}

\begin{table}[h]
\begin{center}
{ \doublerulesep18.0pt \tabcolsep0.01in
\begin{tabular}{|c |c| c | c| }
\hline $m$   & $\sigma_{n+1}/\sigma_n$ &
$\sigma_{n+1}/\hat{\sigma}_n$& $\sigma_1/\hat{\sigma}_n$
 \\
\hline
 $~200~$   & $~7.07\times 10^{-2}~$ &
 $~ 7.07\times 10^{-1}~$ & $~1.00\times 10^1~$
 \\
 \hline
  $~500~$   & $4.47\times 10^{-2}$ &
 $~ 7.07\times 10^{-1}$ & $1.58\times 10^1$
 \\
 \hline
  $~1000~$   & $3.16 \times 10^{-2}$ &
 $~ 7.07\times 10^{-1}$ & $ 2.24\times 10^1$
 \\
 \hline
\end{tabular}}
\caption{Example 1}\label{newtab5}
\end{center}
\end{table}

In Table~\ref{tab5}, we list the results of the relative TLS condition number
$\kappa^{\text{rel}}_g(A,b)$ and its bounds
$$
\underline{\kappa}^{\text{rel}}_2 := \frac{\underline{\kappa}_2 }
{\|x_{TLS}\|}\|[A,\,b]\|_F, \,\, \bar{\kappa}^{\text{rel}}_2 :
= \frac{\bar{\kappa}_2 }{\|x_{TLS}\|}
\|[A,\,b]\|_F, \,\, \bar{\kappa}^{\text{rel}}_1 := \frac{\bar{\kappa}_1
}{\|x_{TLS}\|}\|[A,\,b]\|_F,
$$
and $\kappa^{\rm rel}_{TLS}(A,b)$ and $\bar{\kappa}^{\rm rel}(A,b)$
(see (\ref{golubloanbound}) and (\ref{bgrelbound})),
where $\underline{\kappa}_2$ and $\bar{\kappa}_2$ are defined in
(\ref{TLSimprove}) and $\bar{\kappa}_1$ is defined in (\ref{lubound}).
In Table~\ref{newtab5}, we give some important ratios which have effects on
some of the relative condition numbers listed above.

We can see that the test TLS problems are well conditioned. Both the distance of
$\sigma_{n+1}$ and $\hat{\sigma}_n$ and that of $\sigma_{n+1}$ and ${\sigma}_n$
are not very small so that $\underline{\kappa}^{\text{rel}}_2$,
$\bar{\kappa}^{\text{rel}}_2$ and $\bar{\kappa}^{\text{rel}}_1$ all
estimate $\kappa^{\text{rel}}_g(A,b)$ accurately. Since the three
$\sigma_1/\hat{\sigma}_n$ are considerably bigger than
one, it is known from (\ref{compare}) and the comments after it that
our upper bound $\bar{\kappa}_1^{\rm rel}$ is significantly more accurate than
$\bar{\kappa}^{\rm rel}(A,b)$. Table~\ref{tab5} confirms this. Furthermore,
we see that $\kappa^{\rm rel}_{TLS}(A,b)$ and $\bar{\kappa}^{\rm rel}(A,b)$ are
comparable and also good, but they are not as good as our bounds and overestimate
$\kappa^{\text{rel}}_g(A,b)$ by one to two orders.

{\bf{Example 2.}} In this example, we take the TLS problem from
\cite{KammNagy:1998}. Specifically, a lower $m \times (m - 2 \omega) $
Toeplitz matrix $\bar{T}$ is constructed such that the first column
\begin{equation}\nonumber
t_{i,1} = \left\{
  \begin{array}{ll}
    \frac{1}{\sqrt{2 \pi \beta^2} {\rm{exp}}
    \left[ \frac{(\omega -i+1)^2}{2 \beta^2} \right]} &
    \hbox{\,\,\,$i = 1, 2, \ldots, 2 \omega  +1$,} \\
    0 & \hbox{\,\,\,otherwise,}
  \end{array}
\right.
\end{equation}
and the first row is zero except $t_{1,1}$,
where $\beta =1.25$ and $\omega = 8$.  A Toeplitz matrix $A$ and a right-hand
side vector $b$ are constructed as $A = \bar{T} + E$ and $b = \bar{g} +
e$, where $\bar{g} =[1, \ldots, 1]^T$, $E$ is a random Toeplitz
matrix with the same structure as $\bar{T}$ and $e$ is a random
vector. The entries in $E$ and $e$ are generated randomly from a
normal distribution with mean zero and variance one, and scaled so
that
\begin{equation}\nonumber
\|e\| = \gamma \|\bar{g}\|, \,\,\,\|E\| = \gamma \|\bar{T}\|,
\,\,\gamma = 0.001.
\end{equation}


\begin{table}[h]
\begin{center}
{ \doublerulesep18.0pt \tabcolsep0.01in
\begin{tabular}{|c| c | c| c | c  |c|c|}
\hline $m$   &
$\kappa^{\text{rel}}_g(A,b)$ &
 $\underline{\kappa}^{\text{rel}}_{2}$&
 $\bar{\kappa}^{\text{rel}}_{2} $  & $\bar{\kappa}^{\text{rel}}_{1} $  & $\bar{\kappa}^{\text{rel}}(A,b)$ &
 $\kappa^{\rm rel}_{TLS}(A,b) $
 \\
\hline
 $~100~$   &
  $3.41\times 10^7$   & $3.24\times 10^7$&
 $6.68\times 10^8$ &$1.31\times 10^{11}$& $~ 9.66\times 10^{14}$ &
 $4.70\times 10^{17}$
 \\
 \hline
  $~300~$   &
  $1.62\times 10^8$
 & $1.57\times 10^8$& $4.96\times 10^9$ &$1.10\times 10^{12}$ & $~ 2.24\times 10^{16}$
 &$2.54\times 10^{19}$
 \\
 \hline
  $~500~$   &
  $ 9.50\times 10^7$   & $9.13\times 10^7$
 & $4.31\times 10^9$ &$2.46\times 10^{11}$ & $~ 7.55\times 10^{15}$
 &$1.91\times 10^{19}$
 \\
 \hline
\end{tabular}}
\caption{Example 2}\label{tab4}
\end{center}
\end{table}

\begin{table}[h]
\begin{center}
{ \doublerulesep18.0pt \tabcolsep0.01in
\begin{tabular}{|c |c| c | c| }
\hline $m$   & $\sigma_{n+1}/\sigma_n$ &
$\sigma_{n+1}/\hat{\sigma}_n$& $\sigma_1/\hat{\sigma}_n$
 \\
\hline
 $~100~$   & $~0.964~$ &
 $~ 1-6.08\times 10^{-8}~$ & $~1.04\times 10^4~$
 \\
 \hline
  $~300~$   & $0.984$ &
 $~ 1- 2.03\times 10^{-8}$ & $2.87\times 10^4$
 \\
 \hline
  $~500~$   & $0.993$ &
 $~ 1- 1.37\times 10^{-7}$ & $ 4.34\times 10^4$
 \\
 \hline
\end{tabular}}
\caption{Example 2}\label{newtab4}
\end{center}
\end{table}

In this example, for each test TLS problem, we compute the same quantities as
those in Tables~\ref{tab5}--\ref{newtab5}. The results are reported in
Tables~\ref{tab4}--\ref{newtab4}. As indicated by the $\kappa^{\text{rel}}_g(A,b)$'s,
these TLS problems are all ill conditioned. Their ill conditioning is also
reflected by the fact that $\sigma_{n+1}$ and $\hat{\sigma}_n$ are close.
As estimates of the relative condition number $\kappa^{\text{rel}}_g(A,b)$,
both the lower bound $\underline{\kappa}^{\text{rel}}_{2}$ and the upper bound
$\bar{\kappa}^{\text{rel}}_{2}$ are sharp since $\sigma_{n+1}$ and ${\sigma}_n$
are not so close, but the upper bound $\bar{\kappa}^{\text{rel}}_{1}$ is not
tight any longer and overestimates $\kappa^{\text{rel}}_g(A,b)$ by about
four orders. We see that $\bar{\kappa}^{\text{rel}}_{2}$ improves
$\bar{\kappa}^{\text{rel}}_{1}$ by two orders.
Even though it is not satisfying, $\bar{\kappa}^{\text{rel}}_{1}$
is still much better than $\kappa^{\rm rel}_{TLS}(A,b)$ and
$\bar{\kappa}^{\rm rel}(A,b)$, and the latter two severely overestimate
$\kappa^{\text{rel}}_g(A,b)$ by seven to eight orders and ten to twelve orders,
respectively.

\begin{figure}[h]
\begin{center}
\includegraphics[width=10cm,height=7cm]{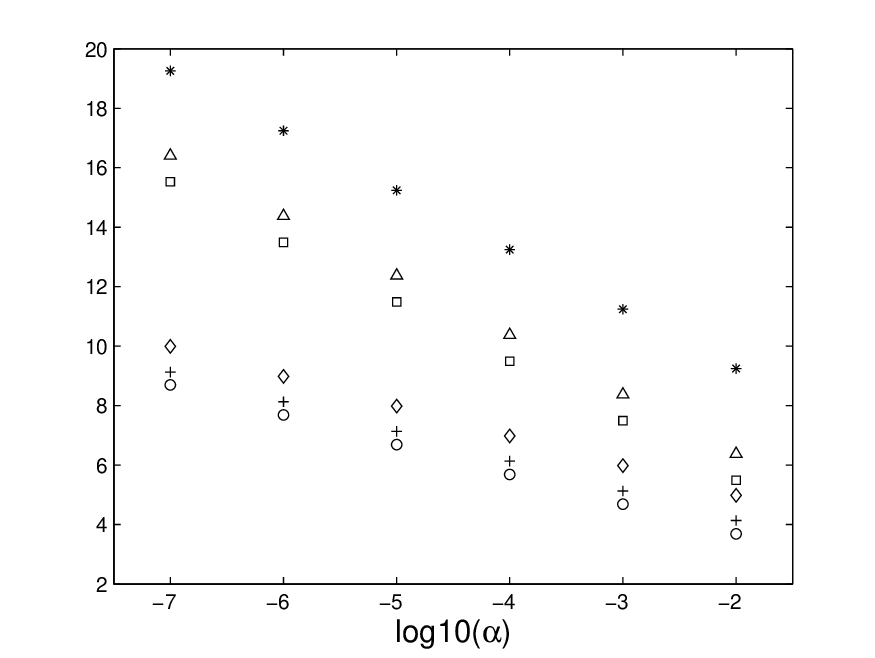}
\end{center}\caption{  $log_{10}(\kappa^{\text{rel}}_g(A,b))$
({\small{${+}$}}),  $log_{10}(\underline{\kappa}^{\text{rel}}_{2})$ ($\circ$),
$log_{10}(\bar{\kappa}^{\text{rel}}_{2})$ ($\diamond$),
$log_{10}(\bar{\kappa}^{\text{rel}}_{1}) $ ({\small{$\Box$}}),
$log_{10}(\bar{\kappa}^{\text{rel}}(A,b)) $ $(\Delta)$
and $log_{10}(\kappa^{\rm rel}_{TLS}(A,b))$
(${*}$) for $(m,n) = (500, 350)$.}\label{figure1}
\end{figure}

{\bf{Example 3.}} Keep in mind that the distance between $\sigma_{n+1}$
and $\hat{\sigma}_n$ can control the conditioning of the TLS problem; see
Theorems \ref{myfirstestimate}--\ref{improveTLSthm} or (\ref{golubloanbound})
and (\ref{bgrelbound}).
In this example, we compare the bounds $\underline{\kappa}^{\text{rel}}_{2}$,
$\bar{\kappa}^{\text{rel}}_{2}$, $\bar{\kappa}^{\text{rel}}_{1} $,
$\kappa^{\rm rel}_{TLS}(A,b)$ and $\bar{\kappa}^{\rm rel}(A,b)$
for various distances between $\sigma_{n+1}$ and $\hat{\sigma}_n$.
On the other hand, keep (\ref{alpha}) in mind. Lemma 4.3 in \cite{GolubVanLoan:1980}
gives
\begin{equation}\label{Golubineqn}
\frac{|\hat{u}^T_n b |}{2 \|x_{TLS}\|} \leq \hat{\sigma}_n - \sigma_{n+1}
\leq \frac{ \|b\|}{ \|x_{TLS}\|},
\end{equation}
which tells us that a small $\alpha$ implies that $\hat{\sigma}_n$
and $\sigma_{n+1}$ are close in some sense. In view of it, for given $m, n$,
we construct $A\in \mathbb{R}^{m \times n}$ and $b \in \mathbb{R}^{n}$ with
different distances between $\sigma_{n+1}$ and $\hat{\sigma}_n$ by
taking different values of $\alpha$.

To do this, we first generate two $n \times n$ random orthogonal matrices
$\bar{U}$ and $\bar{V}$ in the standard normal distribution. Then we
take $\alpha = 10^{-2}, 10^{-3}, \ldots, 10^{-7}$ and run the following function
\begin{eqnarray}\nonumber
&&V = {\text{{{generator}}}}(n,  \alpha, \bar{U}, \bar{V}) \\ \nonumber
&& u = \bar{U}(:,n); \\\nonumber
&& vt = \bar{V}(:,n)';\\\nonumber
&& V_{11} = \bar{U}(:,1:n-1) * \bar{V}(:,1:n-1)' + \alpha * u * vt;\\\nonumber
&& t = {\text{sqrt}}(1 - \alpha * \alpha); \\\nonumber
&& V = \left[~[V_{11}, t * u]; [t * vt, -\alpha]~\right];
\end{eqnarray}
respectively. In such a way, we get six orthogonal matrices
$V \in \mathbb{R}^{(n+1) \times (n+1)}$ with $V(n+1, n+1)=-\alpha$ and
$\alpha=10^{-2}, 10^{-3}, \ldots, 10^{-7}$, respectively.
The idea of construction comes from Proposition~\ref{GolubCSdecomp}.
We generate randomly one $m\times (n + 1)$ matrix $C_1$, and
compute $C_1 =U\Sigma V_1^T$ , the thin SVD of $C_1$. With the matrices
$U$ and $\Sigma$ unchanged, we construct six matrices $C=U\Sigma V^T$ by replacing
$V_1$ by the six orthogonal matrices $V$'s generated above.
Set $A=C(:, 1:n)$, $b = C(:, n+1)$, respectively. Then we get
six different TLS problems.
For each of them, $[A,\, b] = U \Sigma V^T$ is the thin SVD of $[A,\, b]$
and $V(n+1, n+1)=-\alpha$, where $\alpha=10^{-2}, 10^{-3}, \ldots, 10^{-7}$,
respectively.

\begin{figure}[h]
\begin{center}
\includegraphics[width=10cm,height=7cm]{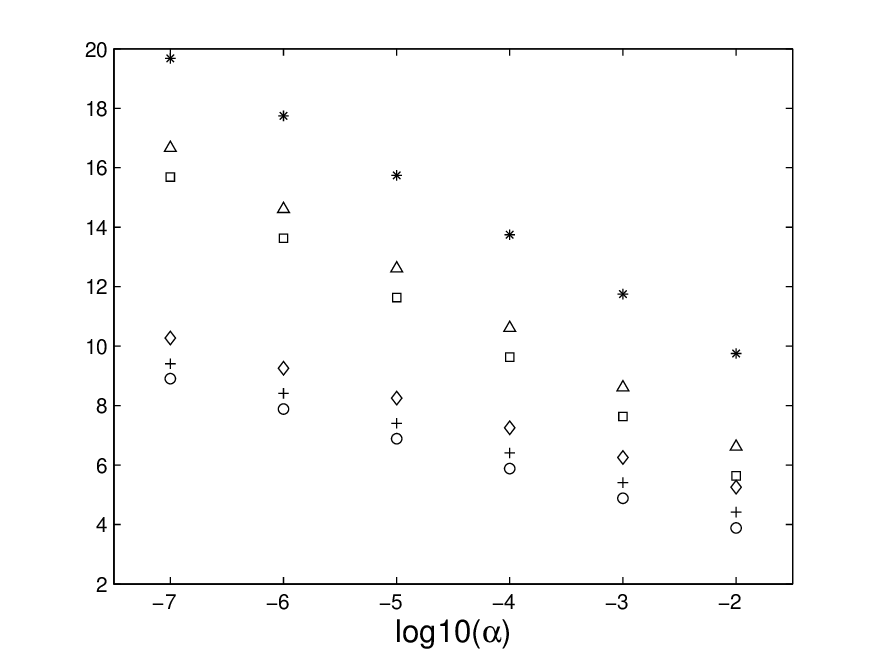}
\end{center}\caption{   $log_{10}(\kappa^{\text{rel}}_g(A,b))$
({\small{$+$}}),
$log_{10}(\underline{\kappa}^{\text{rel}}_{2})$ ($\circ$),
$log_{10}(\bar{\kappa}^{\text{rel}}_{2})$ ($\diamond$),
$log_{10}(\bar{\kappa}^{\text{rel}}_{1}) $ ({\small{$\Box$}}),
$log_{10}(\bar{\kappa}^{\text{rel}}(A,b)) $ $(\Delta)$
and $log_{10}(\kappa^{\rm rel}_{TLS}(A,b))$ (${*}$) for $(m,n) = (1000, 750)$.}
\label{figure2}
\end{figure}

In such a way, with $(m,n) = (500, 350)$ and $(1000, 750)$,
we generate $100$ samples for each $\alpha$, respectively.
For each set of TLS problems with the same $\alpha$,
we compute $\kappa^{\text{rel}}_g(A,b)$, $\underline{\kappa}^{\text{rel}}_{2}$,
$\bar{\kappa}^{\text{rel}}_{2} $, $\bar{\kappa}^{\text{rel}}_{1} $,
$\kappa^{\rm rel}_{TLS}(A,b)$ and $\bar{\kappa}^{\rm rel}(A,b)$.
We plot the ($\log$ scale) averages of these quantities and
the corresponding (log scale) $\alpha$ in Figures~\ref{figure1}--\ref{figure2}.
We also report the averages of $1- \frac{\sigma_{n+1}}{\hat{\sigma}_n}$,
a measure of the distance between  $\sigma_{n+1}$ and $\hat{\sigma}_n$,
in Table~\ref{tab2}. We comment that the averages of
$\frac{\sigma_{n+1}}{\sigma_n}$ for $(m,n) = (500, 350)$ and
$(1000, 750)$ are $0.943$ and $0.958$, respectively.

\begin{table}[h]
\begin{center}
{ \doublerulesep40.0pt \tabcolsep0.01in
\begin{tabular}{c|c|c|c|c|c|c }
\hline
\multicolumn{7}{c}{${(m,n)=(500, 350)}$} 
\\
\cline{1-7}
$\alpha$
 & $~10^{-2}~$ & $~10^{-3}~$& $~10^{-4}~$ &
$~10^{-5}~$ & $~10^{-6}~$ & $~10^{-7}~$
 \\
 \hline
 $1- \frac{\sigma_{n+1}}{\hat{\sigma}_n}$ & $2.91\times 10^{-4}$
 & $2.92\times 10^{-6}$
 & $2.92\times 10^{-8}$& $2.92\times 10^{-10}$ & $2.92\times 10^{-12}$ &
 $2.88\times 10^{-14}$
\\
\hline
\multicolumn{7}{c}{${(m,n)=(1000, 750)}$}
\\
\cline{1-7}
$\alpha$
 & $10^{-2}$ & $10^{-3}$& $10^{-4}$ &
$10^{-5}$ & $10^{-6}$ & $10^{-7}$
 \\
 \hline
 $1- \frac{\sigma_{n+1}}{\hat{\sigma}_n}$ & $3.65\times 10^{-4}$ &
 $3.66\times 10^{-6}$ & $3.66\times 10^{-8}$& $3.66\times 10^{-10}$ &
 $3.66\times 10^{-12}$ & $3.52\times 10^{-14}$
\\
\hline
\end{tabular}}
\caption{Example 3}\label{tab2}
\end{center}
\end{table}

We can see from Figures~\ref{figure1}--\ref{figure2} and Table~\ref{tab2}
that as $\alpha$ decreases, $\hat{\sigma}_n$ and $\sigma_{n+1}$ become closer,
and the TLS problem becomes worse conditioned. $\kappa^{\rm rel}_{TLS}(A,b)$
always severely overestimates $\kappa^{\rm rel}_g(A,b)$. For $\alpha = 10^{-2}$
in which $\hat{\sigma}_n$ and $\sigma_{n+1}$ are not very close,
$\bar{\kappa}^{\rm rel}_1$ is tight and estimate $\kappa^{\rm rel}_g(A,b)$
accurately. For $\alpha=10^{-3}$, $\hat{\sigma}_n$ and $\sigma_{n+1}$ are
closer. In this case, $\bar{\kappa}^{\rm rel}_1$ is no longer tight and estimate
$\kappa^{\rm rel}_g(A,b)$ poorly but it still improves $\kappa^{\rm rel}_{TLS}(A,b)$
and $\bar{\kappa}^{\rm rel}(A,b)$ by about four orders and one order, respectively.
We observe from Figures~\ref{figure1}--\ref{figure2}
that $\bar{\kappa}^{\rm rel}_1$, $\kappa^{\rm rel}_{TLS}(A,b)$ and
$\bar{\kappa}^{\rm rel}(A,b)$ estimate $\kappa^{\rm rel}_g(A,b)$ more poorly as
$\alpha$ decreases. Even so, $\bar{\kappa}^{\rm rel}_1$ is always smaller than
$\kappa^{\rm rel}_{TLS}(A,b)$ and $\bar{\kappa}^{\rm rel}(A,b)$ by about four
orders and one order, respectively.
Remarkably, for all the cases, since ${\sigma}_n$ and $\sigma_{n+1}$ are not so
close, $\underline{\kappa}^{\rm rel}_{2}$ and $\bar{\kappa}^{\rm rel}_{2}$
always estimate $\kappa^{\rm rel}_g(A,b)$ accurately.

\begin{table}[h]
\begin{center}
{ \doublerulesep18.0pt \tabcolsep0.01in
\begin{tabular}{|c |c| c | c| c | c |c |}
\hline $(m,n)$   
& $\kappa^{\text{rel}}_g(A,b)$ &
 $\underline{\kappa}^{\text{rel}}_{2}$&
 $\bar{\kappa}^{\text{rel}}_{2} $  & $\bar{\kappa}^{\text{rel}}_{1} $  & $\bar{\kappa}^{\text{rel}}(A,b)$ &
 $\kappa^{\rm rel}_{TLS}(A,b) $
 \\
\hline
 $~(200,75)~$    & $~1.46\times 10^3~$   & $~5.59\times 10^2~$&
 $~8.58\times 10^3~$ &$~3.19\times 10^{5}~$
 &$~1.01\times 10^{6}~$ &$~2.48\times 10^{8}~$
 \\
 \hline
  $~(500,350)~$   & $9.16\times 10^3$
 & $3.46\times 10^3$& $7.22\times 10^4$ &$1.23\times 10^{6}$
 &$9.45\times 10^{6}$ &$8.16\times 10^{9}$
 \\
 \hline
  $~(1000,75)~$    & $ 3.14\times 10^3$   & $6.19\times 10^2$
 & $2.22\times 10^4$ &$2.57\times 10^{5}$
 &$3.58\times 10^{5}$ &$2.13\times 10^{9}$
 \\
 \hline
  $~(1000,750)~$   & $ 9.73\times 10^4$   & $3.13\times 10^4$
 & $6.41\times 10^5$ &$3.25\times 10^{8}$
 &$3.15\times 10^{9}$ &$1.10\times 10^{12}$
 \\
 \hline
\end{tabular}}
\caption{Example 4}\label{tab6}
\end{center}
\end{table}

{\bf{Example 4.}}
In this example, we generate the entries of $A$ and $b$ as random variables
normally distributed with mean zero and variance one and observe
$\underline{\kappa}^{\text{rel}}_{2}$, $\bar{\kappa}^{\text{rel}}_{2} $,
$\bar{\kappa}^{\text{rel}}_{1} $, $\kappa^{\rm rel}_{TLS}(A,b)$
and $\bar{\kappa}^{\rm rel}(A,b)$.
For each $(m,n)$, we conducted 100 random experiments.
We report the average results of 100 experiments in Table~\ref{tab6}.
We observe that, as estimates of $\kappa^{\text{rel}}_g(A,b)$, both
$\underline{\kappa}^{\text{rel}}_{2}$ and $\bar{\kappa}^{\text{rel}}_{2} $
are tight. The upper bound $\bar{\kappa}^{\text{rel}}_{2} $
improves $\bar{\kappa}^{\text{rel}}_{1} $ by one to two orders
and improves $\kappa^{\rm rel}_{TLS}(A,b)$ and $\bar{\kappa}^{\rm rel}(A,b)$
by about five orders and one to four orders, respectively.
$\bar{\kappa}_1^{\rm rel}$ is always smaller than $\bar{\kappa}^{\rm rel}(A,b)$.
Clearly, the test TLS problems are quite well conditioned, but
$\kappa^{\rm rel}_{TLS}(A,b)$ is a rather poor upper bound and
overestimate $\kappa^{\rm rel}_{g}(A,b)$ too much.

\section{Concluding Remarks}

In the paper, we have studied the SVD-based condition number
theory of the TLS problem. For the TLS condition number, we have established
a new closed formula. Starting with it, we have derived sharp lower and
upper bounds. Importantly and more practically, we have presented
both lower and upper bounds that use only the smallest two singular
values of $A$ and $[A,\ b]$. Numerical experiments have demonstrated
the tightness of our bounds and the improvements of them over the two
upper bounds in \cite{BaboulinGratton:2010,GolubVanLoan:1980}. Throughout
the paper, the considered TLS problem is assumed to satisfy condition
(\ref{PaigeCondition}) and has a unique
TLS solution. It is significant and important to extend the results
presented in the paper to a general generic TLS
problem \cite{HuffelVandewalle:1991,WeiM:1992}
that has non-unique TLS solutions or to the non-generic TLS
problem \cite{HuffelVandewalle:1991}. We will consider these problems
in forthcoming papers. Besides, it might be worthwhile to investigate how to
apply the core problem theory \cite{PaigeStrakos:2006} to study the TLS
condition number.

\bigskip
\bigskip

\leftline {\bf Acknowledgements}

\bigskip

The authors wish to thank the anonymous referees and the editor
Professor Lothar Reichel for their suggestions and comments, which made us
improve the presentation of the paper.

\newcommand{\Gathen}{\relax}\newcommand{\Hoeij}{\relax}


\begin{thebibliography}{10}
\expandafter\ifx\csname url\endcsname\relax
  \def\url#1{\texttt{#1}}\fi
\expandafter\ifx\csname urlprefix\endcsname\relax\def\urlprefix{URL
}\fi


\bibitem{BaboulinGratton:2010}
M. Baboulin, S.~Gratton, {\em A contribution to the conditioning of
the total least squares problem}, SIAM J. Matrix Anal. Appl., 32 (2011)
685--699.

\bibitem{Bjorck:1996}
{\AA}.~Bj\"{o}rck, {\em Numerical Methods For Least Squares
Problems}, SIAM, Philadelphia, PA, 1996.


\bibitem{BjorckHeggernesMatstoms:2000}
{\AA}.~Bj\"{o}rck, P.~Heggernes, P.~Matstoms, {\em Methods for large
scale total least squares problems}, SIAM J. Matrix Anal. Appl., 22
(2000) 413--429.




\bibitem{FierroBunch:1996}
R.~D.~Fierro, J.~R.~Bunch, {\em Perturbation theory for orthogonal
projection methods with applications to least squares and total
least squares}, Linear Algebra Appl., 234 (1996) 71--96.

\bibitem{GohbergKoltracht:1993}
I.~Gohberg, I.~Koltracht, {\em Mixed, componentwise, and structured
condition numbers}, SIAM J. Matrix. Anal. Appl., 14 (1993) 688--704.

\bibitem{GolubVanLoan:1980}
G.~H.~Golub, C.~F.~Van Loan, {\em An analysis of the total least
squares problem}, SIAM J. Numer. Anal., 17 (1980) 883--893.

\bibitem{GolubVanLoan:1996}
G.~H. Golub, C.~F. Van Loan, Matrix Computations, 3rd Edition, Johns
Hopkins University Press, Baltimore, MD, 1996.

\bibitem{Graham:1981}
A.~Graham, Kronecker Products and Matrix Calculus with Application,
Wiley, New York, 1981.

\bibitem{HnetynkovaPlesingerSimaStrakosHuffel:2011}
I.~Hn$\check{e}$tynkov$\acute{a}$, M. Ple$\check{s}$inger,
D. Maria Sima, Z. Strako$\check{s}$, S.~Van Huffel,
{\em The total least squares problem in $AX \approx B$: A new
classification with the relationship to the classical works},
SIAM J. Matrix. Anal. Appl., 32 (2011) 748--770.

\bibitem{Horn:1999}
R.~A.~Horn, C.~R.~Johnson, {\em Matrix Analysis}, Cambridge
University Press, New York, 1985.







\bibitem{jianiu03} Z. Jia, D. Niu, {\em An implicitly restarted refined
bidiagonalization Lanczos method for computing a partial singular value
decomposition}, SIAM J. Matrix Anal. Appl., 25 (2003) 246--265.

\bibitem{jianiu10} Z. Jia, D. Niu, {\em A refined harmonic Lanczos
bidiagonalization method and an implicitly restarted algorithm for
computing the smallest singular triplets of large matrices},
SIAM J. Sci. Comput., 32 (2010) 714--744.

\bibitem{KammNagy:1998}
J.~Kamm, J.~G.~Nagy, {\em A total least squares method for
Toeplitz system of equations}, BIT, 38 (1998) 560--582.

\bibitem{LiJia:2009}
B.~Li, Z.~Jia, {\em Some results on condition numbers of the scaled
total least squares problem}, Linear Algebra Appl., 435 (2011) 674--686.

\bibitem{Liu:1996}
X.~Liu, {\em On the solvability and perturbation analysis for total
least squares problem}, Acta Mathematicae Applicatae Sinica, 19
(1996) 253--262 (in Chinese).

\bibitem{Malyshev:2003}
A.~N.~Malyshev, {\em A unified theory of conditioning for linear
least squares and Tikhonov regularization solutions}, SIAM J.
Matrix. Anal. Appl., 24 (2003) 1186--1196.



\bibitem{PaigeStrakos:2002}
C.~C.~Paige, Z.~Strako$\check{s}$, {\em Scaled total least squares
fundamentals}, Numer. Math., 91 (2002) 117--146.

\bibitem{PaigeStrakos:2006}
C.~C.~Paige, Z.~Strako$\check{s}$, {\em Core problems in linear algebraic systems},
SIAM J. Matrix. Anal. Appl., 27 (2006) 861--875.


%
%

\bibitem{Rice:1966}
J.~R.~Rice, {\em A theory of condition}, SIAM J. Numer. Anal., 3
(1966) 287--310.













\bibitem{HuffelVandewalle:1991}
S.~Van Huffel, J.~Vandewalle, {\em The Total Least Squares Problem:
Computational Aspects and Analysis}, SIAM, Philadelphia, PA, 1991.

\bibitem{WeiM:1992}
M.~Wei, {\em The analysis for the total least squares problem with
more than one solution}, SIAM J. Matrix. Anal. Appl., 13 (1992)
746--763.







\bibitem{ZhouLinWeiQiao:2008}
L.~Zhou, L.~Lin, Y.~Wei, S.~Qiao, {\em Perturbation analysis and
condition numbers of scaled total least squares problems}, Numer.
Algor., 51 (2009) 381--399.

\end{thebibliography}
\end{document}